\documentclass[14pt]{amsart}
\usepackage{lscape, comment} 
\usepackage{enumerate}
\usepackage{xcolor}
\usepackage{amssymb,amsmath,mathtools}
\usepackage{lscape, setspace}
\usepackage{enumerate}
\usepackage{multirow}
\usepackage{booktabs}
\usepackage{latexsym}
\usepackage{amscd,amsfonts}
\usepackage{graphicx}
\usepackage{xcolor}
\usepackage{pdflscape}
\usepackage{longtable}
\usepackage{epsfig}
\usepackage{tikz}
\usepackage{adjustbox}
\usepackage{svg}
\usepackage{float}

\theoremstyle{plain}
\newtheorem{theorem}{Theorem}[section]

\newtheorem{proposition}[theorem]{Proposition}
\newtheorem{lemma}[theorem]{Lemma}
\newtheorem{corollary}[theorem]{Corollary}

\theoremstyle{definition}

\newtheorem{remark}[theorem]{Remark}
\newtheorem{example}[theorem]{Example}

\DeclareMathOperator{\Aut}{Aut}

\DeclareMathOperator{\GL}{GL}

\DeclareMathOperator{\Mat}{Mat}

\newcommand{\F}{\mathbb{F}}
\newcommand{\Z}{\mathbb{Z}}
\newcommand{\N}{\mathbb{N}}

\title[Jordan Superalgebras Type (2,2)]
{The Variety of Jordan Superalgebras of dimension four and even part of dimension two}

\author[Hern\'andez]{ Isabel Hern\'andez}
\address{Secretar\'{\i}a de Ciencia,  Humanidades, Tecnolog\'{\i}a e Innovaci\'on and Centro de Investigaci\'on en Matem\'aticas, Unidad M\'erida}
\email{isabel@cimat.mx }

\author[Martin]{ Mar\'ia Eugenia Martin}
\address{Federal University of Paraná – Curitiba, Brazil}
\email{eugenia@ufpr.br}

\author[Rodrigues]{ Rodrigo Lucas Rodrigues}
\address{Federal University of Cear\'a - Fortaleza, Brazil}
\email{rodrigo@mat.ufc.br}

\thanks{The first author was
supported by Grant CONAHCYT A1-S-45886.}
\begin{document}

\begin{abstract}
We describe the variety of Jordan superalgebras of dimension $4$ whose even part is a Jordan algebra of dimension $2$ over an algebraically closed field $\mathbb{F}$ of characteristic $0$. We prove that the variety has $25$ irreducible components, $24$ of them correspond to the Zariski closure of the $GL_2(\mathbb{F})\times GL_2(\mathbb{F})$-orbits of rigid superalgebras and the other one is the Zariski closure of an union of orbits of an infinite family of superalgebras, none of them being rigid. Furthermore, it is known that the question of the existence of a rigid Jordan algebra whose second cohomology group does not vanish is still an open problem. We solve this problem in the context of superalgebras, showing a four-dimensional rigid Jordan superalgebra whose second cohomology group does not vanish. 
\end{abstract}
\maketitle

\section{Introduction}
Non-necessarily associative algebras have been extensively studied. Among the main classes are the Lie and Jordan algebras. It is essential to know examples of a given algebraic system, particularly in low dimensions. 

In relation to finite-dimensional algebras over an arbitrary field, it is natural to study their algebraic classification, that is, to determine all that are not isomorphic. It is also of interest to determine the deformations between low-dimensional algebras in a variety defined by polynomial identities, in particular finding the rigid ones. The irreducible components within the regarded variety related to Zariski topology are constituted by the closure of their orbits under the action of the linear group.

We now briefly describe the progress of the geometric classification of Jordan algebras and superalgebras. In $2004$, I. Kashuba and I. Shestakov \cite{teseIryna} investigated the variety ${\rm Jor}_3$ of Jordan algebras of dimension $3$, for $\F$ algebraically closed of characteristic $\neq 2, 3$, listing all $\GL_3(\F)$-orbits of ${\rm Jor}_3$ and establishing its $5$ irreducible components. In $2006$, I. Kashuba \cite{iryna} studied the variety ${\rm Jor}_n$ of Jordan unitary algebras of dimension $n$, for an algebraically closed field $\F$ with $\textrm{char\,} \F \neq 2$, as well as infinitesimal deformations of Jordan algebras, establishing the list of $\GL_n(\F)$-orbits on ${\rm Jor}_n$, $n = 4, 5$, under the ``change of basis'' action, finding the number 
of irreducible components as being $3$ and $6$ respectively, and a list of rigid algebras was included. In $2014$, I. Kashuba and M. E. Martin \cite{kashubamartin} studied the variety ${\rm Jor}_4$ of four-dimensional Jordan algebras, for an algebraically closed field $\F$ of characteristic $\neq 2$, describing its irreducible components and proving that ${\rm Jor}_4$ is the union of Zariski closures of the orbits of 10 rigid algebras.
In $2016$, I. Kashuba and M. E. Martin \cite{MartinJor3} studied the variety ${\rm Jor}_3$ of three-dimensional Jordan algebras over the field of real numbers, establishing the list of $26$ non-isomorphic Jordan algebras and describing the irreducible components of ${\rm Jor}_3$, proving that it is the union of Zariski closures of the orbits of eight rigid algebras.
In $2018$, I. Kashuba and M. E. Martin \cite{MartinJorN5} also investigated the variety ${\rm JorN}_5$ of five-dimensional nilpotent Jordan algebras structures over an algebraically closed field, showing that ${\rm JorN}_5$ is the union of five irreducible components, four of them correspond to the Zariski closure of the $\GL_5(\F)$-orbits of four rigid algebras and the other one is the Zariski closure of an union of orbits of an infinite family of algebras, none of them being rigid. 

In the $\mathbb{Z}_2$-graded context, in $2019$, M. A. Alvarez et al. \cite{Geom_SuperJAdim3} described degenerations of three-dimensional Jordan superalgebras over the field of complex numbers $\mathbb{C}$. In particular, they describe all irreducible components in the corresponding varieties. Recently, we \cite{GeomSuperJor4-13-31} described the variety of Jordan superalgebras of dimension $4$ whose even part is a Jordan algebra of dimension $1$ or $3$, proving that the variety is the union of $11$ and $21$ rigid superalgebras, respectively, and in both cases the irreducible components of the varieties were also described. 

In this article, we deal with the geometric classification of Jordan superalgebras of dimension four with even and odd parts of dimension 2 over an algebraically closed field $\mathbb{F}$ of characteristic $0$. In Section 2, we set up notation and terminology, we indicated some useful invariants to guarantee the non-existence of a deformation between two superalgebras, and we introduce the notion of Jordan superalgebras. Section 3 provides a necessary condition for a superalgebra to be rigid and some examples are given. In particular, a rigid Jordan superalgebra whose second cohomology group does not vanish is constructed. In Section 4 we describe deformations between Jordan superalgebras of type (2,2) and characterize the irreducible components of such variety. As a direct consequence, we achieve the geometric classification of associative and supercommutative superalgebras of dimension four and type $(2,2)$, resulting in 6 irreducible components. Additionally, we determine that the subvariety of nilpotent Jordan superalgebras of dimension four and type $(2,2)$ has 3 irreducible components.


\section{Preliminaries}
Let $\mathbb{F}$ be any field of characteristic different from two. A superalgebra $A$ is just a $\Z_2$-graded algebra over $\mathbb{F}$, that is, $A$ is decomposed into a direct sum of subspaces $A=A_0 \oplus A_1$ such that $A_iA_j\subseteq A_{i+j}$, for $i, j \in \Z_2$. Notice that $A_0$ is a subalgebra and $A_1$ is a $A_{0}$-bimodule. The elements in $A_0\setminus\{ 0 \}$ (respectively, in $ A_1 \setminus \{0\}$) are called even (respectively, odd) and the elements in $A_i\setminus\{ 0 \}$, $i \in \Z_2$, are said to be homogeneous and of degree $i$. 
The degree of a homogeneous element $x$ is denoted by $|x|=i$. The pair $(\dim_\F (A_0), \dim_\F (A_1))$ is said to be the type of $A$. 

The morphisms in this category are defined as follows. Let $A$ and $A^\prime$ be superalgebras. A linear map $\Phi\colon A \to A^\prime$ is a morphism of superalgebras if $\Phi$ is even (i.e., $\Phi(A _i) \subset A^\prime_i $, for $ i \in \Z_2$) and $\Phi (xy)= \Phi(x) \Phi(y)$, for all $x,y \in A$.

We now introduce the concept of the Grassmann envelope. Let $\mathcal{G}= {\rm alg}\{1, e_i | 1\leq i, \; e_ie_j=-e_je_i \}$ be the Grassmann algebra. Then, $\mathcal{G}=\mathcal{G}_0 \oplus \mathcal{G}_1$ is an associative $\Z_2$-graded algebra, where $\mathcal{G}_0$ (respectively, $\mathcal{G}_1$) are subspaces spanned by the products of even (respectively, odd) length. In addition, $\mathcal{G}$ is supercommutative, that is, $uv=(-1)^{|u||v|}vu$, for all homogeneous elements $u, v \in \mathcal{G}$. The {\it Grassmann envelope} of a superalgebra $A = A_{0} \oplus A_{1}$ is the algebra 
$$
\mathcal{G}(A)= (\mathcal{G}_0 \otimes A_0) \oplus (\mathcal{G}_1 \otimes A_1),
$$ 
where the multiplication is defined by $(x\otimes u)(y\otimes v):= xy\otimes uv$, for all $x\otimes u \in \mathcal{G}_i \otimes A_i$, $y\otimes v \in \mathcal{G}_j \otimes A_j$, and $i,j \in \Z_2$.

Let $\mathcal{V}$ be a variety of algebras over $\F$ defined by a set of multilinear identities $\{ P_\lambda\}_{\lambda\in \Lambda}$. A superalgebra $A=A_0\oplus A_1$ is called a $\mathcal{V}$-superalgebra if its Grassmann envelope $\mathcal{G}(A)$ lies in $\mathcal{V}$. In particular, if $A =A_0\oplus A_1$ is a $\mathcal{V}$-superalgebra, then the algebra $A_0$ lies in $\mathcal{V}$ and $A_1$ is a $A_0$-bimodule in the class $\mathcal{V}$. 

Let $\{P_\lambda\}_{\lambda\in \Lambda}$ be the set of multilinear identities that defines the variety $\mathcal{V}$. In general, the corresponding set of superidentities $\{ P^s_ \lambda\}_{\lambda\in \Lambda}$ which defines the variety of $\mathcal{V}$-superalgebras can be obtained following the {\it Kaplansky Rule:} If two homogeneous adjacent elements $x,y$ are exchanged, then the corresponding term is multiplied by $(-1)^{|x||y| }$. 

\begin{remark} In some varieties of algebras, we impose some restrictions over the characteristic of the field $\F$ to obtain the multilinear identities.
One of the most common conditions is
${\rm char} (\F)\neq 2$, but in some cases we require others, as for example ${\rm char}(\F)\neq 2,3,5$ in the class of commutative power-associative superalgebras.
\end{remark}

The set of all $\mathcal{V}$-superalgebras of type $(m,n)$ defines an affine variety in $\F^{m^3+3mn^2}$ and it will be denoted by $\mathcal{VS}^{(m, n)}$. In fact, let $V=V_0 \oplus V_1$ be a $\Z_2$-graded vector space with a fixed homogeneous basis
$\left\{e_1,\dots ,e_m, f_1,\dots ,f_n\right\}$. Our goal is to provide a superalgebra structure for $V$. To this end, the possible multiplication tables are given by
\begin{eqnarray*}
e_i e_j = \sum_{k=1}^m \alpha_{ij}^k e_k,&\quad& f_i f_j = \sum_{k=1}^m \gamma_{ij}^k e_k,\\
e_i f_j = \sum_{k=1}^n  \beta_{ij}^kf_k, &\quad& 
f_i e_j = \sum_{k=1}^n  {\beta^\prime}_{ij}^k f_k. 
\end{eqnarray*}
where $\alpha_{ij}^k, \beta_{ij}^k, {\beta^\prime}_{ij}^k, \gamma_{ij}^k \in \F$, for all $i, j$ and $k$. Notice that every set of structure constants must satisfy the polynomial identities $P^s_\lambda$, for $\lambda\in \Lambda$. 
It follows that each point $(\alpha_{ij}^k, \beta_{ij}^k, {\beta^\prime}_{ij}^k, \gamma_{ij}^k)\in \mathcal{VS}^{(m, n)}$ represents, in the fixed basis, a $\mathcal{V}$-superalgebra $A$ over $\F$ of type $(m,n)$. 

Since superalgebras morphisms are even maps, there is a natural ``change of basis'' action of the group $G= \GL_m(\F) \times \GL_n(\F)$ on $\mathcal{VS}^{(m, n)}$ which gives a one-to-one correspondence between $G$-orbits $A^G$ on $\mathcal{VS}^{(m, n)}$ and the isomorphism classes of $\mathcal{V}$-superalgebras of type $(m,n)$. Let $\overline{A^G}$ be the Zariski closure of $A^G$. 
The superalgebra $A \in \mathcal{VS}^{(m, n)}$ is said to be a {\it deformation} of $B \in \mathcal{VS}^{(m, n)}$ if
$B^G \subseteq \overline{A^G}$. We indicate it by $A \to B$. The superalgebra $A\in \mathcal{VS}^{(m, n)}$ is called {\it (geometrically) rigid} if the orbit $A^G$ is Zariski-open set in $\mathcal{VS}^{(m, n)}$. Hence, if $A$ is rigid, then $\overline{A^G}$ is an irreducible component of the variety $\mathcal{VS}^{(m, n)}$ and 
any deformation $C$ of $A$ satisfies $C\simeq A$.

The one-parameter family of deformations technique is used in order to describe the deformations on $\mathcal{V}$-superalgebras of type $(m, n)$, and it consists of the following. 
Let $A \in \mathcal{VS}^{(m, n)}$, denote the Laurent polynomials in the variable $t$ by $\F (t)$, and take $g(t)\in \Mat_m(\F(t)) \times \Mat_n(\F(t))$. For any $0 \neq t \in \F$, assume that $g(t) \in G$. If $A_t=(\alpha_{ij}^k(t), \beta_{ij}^k(t), {\beta^\prime}_{ij}^k(t),  \gamma_{ij}^k(t))$ is the $\mathcal{V}$-superalgebra resulting from the application of the change of basis $g(t)$ to $A$, then $A$ is a deformation of the $\mathcal{V}$-superalgebra 
$B=(\alpha_{ij}^k(0), \beta_{ij}^k(0), {\beta^\prime}_{ij}^k(0), \gamma_{ij}^k(0))$. In particular, a $\mathcal{V}$-superalgebra 
$A$ is {\it rigid} if any $g(t)$ satisfying the above conditions defines the $\mathcal{V}$-superalgebra $A_t$ isomorphic to $A$ for every $t \in \F$. In other words, the following result holds. 

\begin{lemma}\label{Lema:curva-Mazzola} \cite{mazzola} A curve $g(t)$ in $\mathcal{VS}^{(m,n)}$ which generically lies in a subvariety $U$ and which cuts $A^G$ in special points implies that $A$ belongs to $\overline{U}$ (and conversely).
\end{lemma}

The following results will be useful to guarantee the non-existence of a deformation between two superalgebras of $\mathcal{VS}^{(m,n)}$. They were demonstrated in \cite{DegLieSuperalg} and \cite{Geom_SuperJAdim3} for the case of Lie and Jordan superalgebras, respectively, but the proofs can be easily adapted to any other variety.

\begin{lemma} \label{Lema:nondeformation}
Let $A, B \in \mathcal{VS}^{(m,n)}$ be such that $A$ has structure constants $(\alpha_{ij}^k, \beta_{ij}^k, {\beta^\prime}_{ij}^k, \gamma_{ij}^k)$ and $A \to B$. The following conditions hold: 
\begin{enumerate}[(i)]
\item \label{dimAut}$\dim(\Aut(A)) < \dim(\Aut(B))$. 
\item \label{dimJ_1^2} $ \dim( A^r)_i \geq \dim (B^r)_i$, where $i\in \mathbb{Z}_2$, $A^1 = A$,
$A^r = A^{r-1}A + \cdots + A A^{r-1}$, for all $r\in \N$, and $A^r = (A^r)_0 \oplus (A^r)_1$.
\item \label{partepar}$A_0 \to B_0$.
\item $ab(A) \to ab(B)$, where $ab(A) \in \mathcal{VS}^{(m,n)}$ has structure constants $(0, 0, 0, \gamma_{ij}^k)$.
\item \label{esquecimento}$\mathcal{F}(A) \to \mathcal{F}(B)$, where $\mathcal{F}(A) \in \mathcal{VS}^{(m,n)}$ has structure constants $(\alpha_{ij}^k, \beta_{ij}^k, 
{\beta^\prime}_{ij}^k,0)$.
\item If $A$ is associative, then $B$ is also associative. Furthermore, if $A$ satisfies a polynomial identity, 
then $B$ satisfies the same polynomial identity.
\end{enumerate} 
\label{lemma:general_invariants}
\end{lemma}

\begin{remark} \label{Remark:como algebras} Suppose that $A, B \in\mathcal{VS}^{(m,n)}$ have structure of $\mathcal{V}$-algebras. If $A \not \to B$ as algebras, then $A \not \to B$ as superalgebras.
\end{remark}

In this article, we consider the variety $\mathcal{V}$ of \textit{Jordan algebras} over an algebraically closed field $\F$ of $\textrm{char\,}\F = 0$, i.e., the variety of algebras defined by the multilinear identities $xy=yx$ and
\[
(wx)(yz)+(wy)(xz)+(wz)(xy)-x(w(yz)) - y(w(xz)) - z(w(xy)) =0.
\]
The last one is known as Jordan identity. We will work with the corresponding variety of $\mathcal{V}$-superalgebras, called \textit{Jordan superalgebras}, whose set of superidentities, obtained through Kaplansky's Rule, are: the supercommutativity $xy=(-1)^{|x||y|} yx$, and the Jordan superidentity:
\begin{eqnarray*}
(wx)(yz)+(-1)^{|x||y|}(wy)(xz)+(-1)^{(|x|+ |y|)|z| }(wz)(xy)& \\
-(-1)^{ |w| |x|} x (w(yz ))- (-1)^{| y| (|w| +|x| )} y(w(xz)) - (-1)^{ |z|(|w|+|x|+|y|)} z(w(xy))=0,
\end{eqnarray*}
where $x,y,z,w$ are homogeneous elements. 
The variety of all Jordan superalgebras of type $(m,n)$ will be denoted by $\mathcal{JS}^{(m,n)}$.

In the next section we describe a powerful machinery to find rigid superalgebras. This technique is valid for other varieties, such as associative and Lie superalgebras. In order to be precise, we introduce it only for Jordan superalgebras.

\section{Second cohomology group of finite-dimensional Jordan Superalgebras}
The most known sufficient condition for a superalgebra to be rigid is given in terms of its cohomology group. 
We say that the second cohomology group $H^2(\mathcal{J},\mathcal{J})$ of a Jordan superalgebra $J$ with coefficients in itself 
\textit{vanishes} if for every bilinear map $h\colon \mathcal{J}\times \mathcal{J}\to \mathcal{J}$ satisfying 
\begin{equation}\label{eq:super-simetria}
h(a,b)=(-1)^{|a||b|}h(b,a)
\end{equation} and 
\begin{align}\label{eq:2-coclico}
\nonumber F(a,b,c,d)+(-1)^{|b|(|c|+|d|)+|c||d|}F(a,d,c,b)+(-1)^{|a|(|b|+|c|+|d|)+|c||d|}F(b,d,c,a)&\\
=G(a,b,c,d)+(-1)^{|b||c|}G(a,c,b,d)+(-1)^{|d|(|c|+|b|)}G(a,d,b,c)&
\end{align}
for all homogeneous elements $a,b,c,d$, 
where 
\begin{equation*}
\begin{array}{c}
F(a,b,c,d):=h((ab)c,d)+h(ab,c)d+(h(a,b)c)d,\\
G(a,b,c,d):=h(ab,cd)+h(a,b)(cd)+(ab)h(c,d),
\end{array}
\end{equation*} there exists a linear map $\mu\colon \mathcal{J} \to \mathcal{J} $ such that
\begin{equation}\label{hmu}
h(a,b)=a \mu(b)+\mu(a) b-\mu(ab).
\end{equation}
For the precise definition and properties of this group for Jordan superalgebras, we refer the reader to \cite{FaberCohomologia} and \cite{FaberCohomologia2}. 

\begin{proposition}\label{cohomology} Let
$ \mathcal{J} \in \mathcal{JS}^{(n,m)}$ be a Jordan superalgebra. If $H^2(\mathcal{J},\mathcal{J})=0$, then $\mathcal{J}$ is a rigid superalgebra. 
\end{proposition}

The last implication applies equally well to any category of algebras or superalgebras defined by superidentities  
with appropriate modifications. In fact, this result was originally obtained in \cite{gerstenhaber} for associative and Lie algebras; for the case of Jordan algebras, see \cite{iryna}.The proof for Jordan superalgebras is analogous.

\begin{example}\label{ex:(2,2)_3_rigid} 
Consider $\mathcal{J}=(2,2)_3\in \mathcal{JS}^{(2,2)}$ whose multiplication table is given in Table \ref{table:JSA_(2,2)} and let $h:\mathcal{J}\times \mathcal{J}\to \mathcal{J}$
be a bilinear map satisfying the supersymmetry \eqref{eq:super-simetria} and the condition \eqref{eq:2-coclico}, then 
\begin{gather*}
h(e_{1},e_{1})=\alpha e_1 + \beta e_2,\quad 
h(e_{1},e_{2})=\gamma e_1 - \beta e_2,\quad 
h(e_{1},f_{1})=\delta e_2 + \frac{1}{2} \alpha f_1,\\ 
h(e_{1},f_{2})=\epsilon e_2 + \frac{1}{2} \alpha f_2,\quad 
h(e_{2},e_{2})=-\gamma e_1 + \lambda e_2,\quad 
h(e_{2},f_{1})=-2 \delta e_2 + \frac{1}{2} \gamma f_1, \\
h(e_{2},f_{2})=-2 \epsilon e_2 + \frac{1}{2} \gamma f_2,\quad h(f_{1},f_{2})= \theta e_1 + \beta e_2,
\end{gather*}
for any $\alpha,\beta,\gamma, \delta, \epsilon, \lambda, \theta\in\mathbb{F}$. Note that by supersymmetry, it follows that $h(f_i,f_i)=0$, for $i=1,2$. Now, define a linear map $\mu\colon J\to J$
by
\begin{gather*}
\mu(e_{1})=\alpha e_{1}-\beta e_{2}, \quad
\mu(e_{2})=\gamma e_{1}+\lambda e_{2}, \\
\mu(f_{1})=-2\delta e_{2}, \quad
\mu(f_{2})=-2\epsilon e_{2}+(\theta+\alpha)f_2.
\end{gather*}
Thus, we have that $\eqref{hmu}$ holds and then $H^{2}(\mathcal{J},\mathcal{J})=0$,
which implies that $(2,2)_3$ is a rigid Jordan superalgebra.
\end{example}

\begin{example} \label{ex:(2,2)_5_rigid}
Analogously, for $ \mathcal{J}=(2,2)_5$ (see Table \ref{table:JSA_(2,2)}) we will calculate the second cohomology group. Let $h: \mathcal{J}\times \mathcal{J} \to \mathcal{J} $
be a bilinear map as before, then 
\begin{gather*}
h(e_{1},e_{1})=\alpha e_1 - \beta e_2,\quad 
h(e_{1},e_{2})=\gamma e_1 + \beta e_2,\quad 
h(e_{1},f_{1})=\delta e_2 + \alpha f_1,\\ 
h(e_{1},f_{2})=\epsilon e_2 + \alpha f_2,\quad 
h(e_{2},e_{2})=-\gamma e_1 + \lambda e_2,\quad 
h(e_{2},f_{1})=- \delta e_2 + \gamma f_1, \\
h(e_{2},f_{2})=- \epsilon e_2 + \gamma f_2,\quad
h(f_{1},f_{2})= \theta e_1 - \beta e_2
\end{gather*}
for any $\alpha,\beta,\gamma, \delta, \epsilon, \lambda, \theta\in\mathbb{F}$. Define a linear map $\mu\colon \mathcal{J} \to \mathcal{J}$
by
\begin{gather*}
\mu(e_{1})=\alpha e_{1}+\beta e_{2}, \quad
\mu(e_{2})=\gamma e_{1}+\lambda e_{2}, \\
\mu(f_{1})=-\delta e_{2}, \quad
\mu(f_{2})=-\epsilon e_{2}+(\theta+\alpha)f_2. 
\end{gather*}
Then $\eqref{hmu}$ holds and hence $H^{2}( \mathcal{J} , \mathcal{J})=0$. Therefore, $(2,2)_5$ is a rigid Jordan superalgebra.
\end{example}

 In positive characteristic, it has been shown that the reverse implication of Proposition \ref{cohomology} is false for finite-dimensional associative algebras. M. Gerstenhaber and S. Schack \cite{gerstenBockstein} constructed rigid high-dimensional associative algebras $A$ over a field of positive characteristic such that $H^2(A, A)\neq 0$. Furthermore, Richardson \cite{Richardson} showed that there exist complex Lie algebras in every even dimension greater than $16$ that are rigid but the second cohomology group does not vanish. The question of the existence of a rigid Jordan algebra $\mathcal{J}$ satisfying $H^2 (\mathcal{J}, \mathcal{J}) \neq 0$ is still an open problem. 
In what follows, we will show that the converse does not hold for Jordan \textit{superalgebras}. We will show an example of a Jordan superalgebra that is rigid (as shown in Theorem \ref{Thm: principal}) but whose second cohomology group does not vanish.
\begin{example}Let $ \mathcal{J}=(2,2)_8$ (see Table \ref{table:JSA_(2,2)}) and consider a bilinear map $h: \mathcal{J}\times \mathcal{J} \to \mathcal{J} $ satisfying \eqref{eq:super-simetria} and \eqref{eq:2-coclico}, then 
\begin{gather*}
h(e_{1},e_{1})=\alpha e_1 + \beta e_2 + \omega f_2,\quad 
h(e_{1},e_{2})=-\eta e_1 - \beta e_2 + \sigma f_1, \\
h(e_{1},f_{1})= \delta e_2 + \frac{1}{2} \alpha f_1 + \rho f_2, \quad
h(e_{1},f_{2})= \pi f_1,\\
h(e_{2},e_{2})=\eta e_1 + \lambda e_2 - 2 \sigma f_1 + \tau f_2,\quad 
h(e_{2},f_{1})=-2 \delta e_2 - \frac{1}{2} \eta f_1, \\
h(e_{2},f_{2})=\nu e_2,\quad
h(f_{1},f_{2})=\kappa f_1
\end{gather*}
for any $\alpha,\beta,\delta, \lambda, \omega, \eta, \sigma, \tau, \rho, \pi, \nu, \kappa \in\mathbb{F}$. Suppose that there exists a linear map $\mu\colon \mathcal{J} \to \mathcal{J}$ where $\mu(e_i)=\sum_{j=1}^4u_{ij}e_j$, with $e_3=f_1$ and $e_4=f_2$, and such that $\eqref{hmu}$ holds.
Taking $a=f_1$ and $b=f_2$ in $\eqref{hmu}$, we get 

\begin{equation*}
\left(\kappa - \frac{1}{2} u_{41}\right)f_1 = 0,
\end{equation*}
and on the other hand, setting $a = f_2$ and $b = f_1$ in $\eqref{hmu}$, we obtain
\begin{equation*}
\left(- \kappa -\frac{1}{2} u_{41}\right)f_1=0,
\end{equation*}
which implies $\kappa=0$. This means that if $\kappa \neq 0$, then there is no such $\mu$, and therefore $\mathcal{J}$ is a rigid Jordan superalgebra such that $H^2(\mathcal{J},\mathcal{J})$ does not vanish. 
\end{example}


Analogously, we say that the \textit{even part} of the second cohomology group $H^2(\mathcal{J},\mathcal{J})$ of a Jordan superalgebra $\mathcal{J}$ with coefficients in itself
\textit{vanishes} (i.e., $(H^2(\mathcal{J},\mathcal{J}))_0=0$) if for every \textit{even} bilinear map $h_0:\mathcal{J}\times \mathcal{J}\to \mathcal{J}$ (this means that $h_0(e_i,e_j),h_0(f_i,f_j)\in J_0$ and $h_0(e_i,f_j),h_0(f_i,e_j)\in \mathcal{J}_1$ for $i,j=1,2$) satisfying \eqref{eq:super-simetria} and \eqref{eq:2-coclico} there exists an \textit{even} linear map $\mu_0: \mathcal{J} \to \mathcal{J} $ (i.e., $\mu_0(\mathcal{J}_i)\subseteq \mathcal{J}_i$ for $i=1,2$) such that
\begin{equation}\label{even-hmu}
 h_0(a,b)=a \mu_0(b)+\mu_0(a) b-\mu_0(ab).
\end{equation}

In \cite{DegLieSuperalg}, the authors proved that it suffices for the even part of the second cohomology group of a Lie superalgebra to vanish in order to conclude that it is rigid. The proof for Jordan superalgebras is analogous. 

\begin{proposition}\label{even-cohomology} If
$ \mathcal{J} \in \mathcal{JS}^{(n,m)}$ and 
the even part of the second cohomology group $(H^2(\mathcal{J},\mathcal{J}))_0$ vanishes, then $\mathcal{J}$ is a rigid Jordan superalgebra. \end{proposition}

We leave open the natural question that arises from Proposition \ref{even-cohomology}, namely whether there exists a rigid Jordan superalgebra $\mathcal{J}$ such that $(H^2 (\mathcal{J}, \mathcal{J}))_0 \neq 0$. 

\begin{example}\label{ex:(2,2)_1_rigid} 
Consider $\mathcal{J}=(2,2)_1\in \mathcal{JS}^{(2,2)}$ (see Table \ref{table:JSA_(2,2)})  and let $h_0\colon \mathcal{J}\times \mathcal{J}\to \mathcal{J}$
be an even bilinear map satisfying the above conditions, then 
\begin{gather*}
h_0(e_{1},e_{1})=\alpha e_1 + \beta e_2,\quad 
h_0(e_{1},e_{2})=\gamma e_1 - \beta e_2,\quad
h_0(e_{1},f_{1})=0,\\ 
h_0(e_{1},f_{2})=0,\quad 
h_0(e_{2},e_{2})=-\gamma e_1 + \delta e_2,\quad 
h_0(e_{2},f_{1})=0, \\
h_0(e_{2},f_{2})=0,\quad h_0(f_{1},f_{2})=0
\end{gather*}
for any $\alpha,\beta,\gamma, \delta \in\mathbb{F}$. Define an even linear map $\mu_0 \colon \mathcal{J}\to \mathcal{J}$ by 
\begin{gather*}
\mu_0(e_{1})=\alpha e_{1}-\beta e_{2}, \quad
\mu_0(e_{2})=\gamma e_{1}+\delta e_{2}, \\
\mu_0(f_{1})=0, \quad
\mu_0(f_{2})=0,
\end{gather*}
thus we have that $\eqref{even-hmu}$ holds and then $(H^{2}(\mathcal{J},\mathcal{J}))_0=0$
which implies $(2,2)_1$ is a rigid Jordan superalgebra.
\end{example}

\section{Jordan Superalgebras of Type $(2,2)$}
In this section, we investigate 
the variety $\mathcal{JS}^{(2,2)}$ over an algebraically closed field of characteristic zero. I. Hernández et al. \cite{SuperPowerAssoc4} have provided a concrete list of non-isomorphic commutative power-associative superalgebras up to dimension $4$ over an algebraically closed field of characteristic prime to $30$. As a consequence of this classification, we see that there exist, up to isomorphism, $72$ Jordan superalgebras of type $(2,2)$ and an one parameter family over an algebraically closed field of characteristic zero (see \cite{SuperPowerAssoc4}, Tables 6, 8, 9, 10, 12, and 14). 
Table \ref{table:JSA_(2,2)} gives representatives for isomorphism classes and some additional useful information, namely the dimension of the automorphism group of each superalgebra, and we
indicate by ``A'' if the superalgebra is associative and ``NA'' otherwise. 
The superscript ``N'' in $(2,2)_i^N$ indicates that the superalgebra is nilpotent.

\begin{longtable}[H]{llcc}
\caption{\label{table:JSA_(2,2)}Jordan superalgebras of type $(2,2)$.}\\
Label & Multiplication table & $\dim(\Aut(\mathcal{J})) $ & \\
\hline
\endhead
$(2,2)_1$: & $\;\;e^2_1=e_1$, $\;\;e^2_2=e_2$. & $4$ & A  \\
$(2,2)_2$: & $\;\;e^2_1=e_1$, $\;\;e^2_2=e_2$, $\;\;e_1 f_1= \frac{1}{2} f_1$, $\;\;e_1 f_2= \frac{1}{2} f_2$. & $4$ & NA \\
$(2,2)_3$: & $\;\;e^2_1=e_1$, $\;\;e^2_2=e_2$, $\;\;e_1 f_1= \frac{1}{2} f_1$, $\;\;e_1 f_2= \frac{1}{2} f_2$, 
$\;\;f_1f_2= e_1$. & $3$ & NA \\
$(2,2)_4$: & $\;\;e^2_1=e_1$, $\;\;e^2_2=e_2$, $\;\;e_1 f_1= f_1$, $\;\;e_1 f_2= f_2$. & $4$& A \\
$(2,2)_5$: & $\;\;e^2_1=e_1$, $\;\;e^2_2=e_2$, $\;\;e_1 f_1= f_1$, $\;\;e_1 f_2= f_2$, 
$\;\;f_1f_2= e_1$. & $3$ & NA  \\
$(2,2)_6$: & $\;\;e^2_1=e_1$, $\;\;e^2_2=e_2$, $\;\;e_1 f_1= \frac{1}{2} f_1$, $\;\;e_1 f_2= \frac{1}{2} f_2$, $\;\;e_2 f_1= \frac{1}{2} f_1$, & $4$ & NA  \\
& $\;\;e_2 f_2= \frac{1}{2} f_2$. &&\\
$(2,2)_7$: & $\;\;e^2_1=e_1$, $\;\;e^2_2=e_2$, $\;\;e_1 f_1= \frac{1}{2} f_1$, $\;\;e_1 f_2= \frac{1}{2} f_2$, $\;\;e_2 f_1= \frac{1}{2} f_1$, & $3$ & NA  \\
& $\;\;e_2 f_2= \frac{1}{2} f_2$, $\;\;f_1f_2= e_2$. &&\\
$D_\gamma$: & $\;\;e^2_1=e_1$, $\;\;e^2_2=e_2$, $\;\;e_1 f_1= \frac{1}{2} f_1$, $\;\;e_1 f_2= \frac{1}{2} f_2$, $\;\;e_2 f_1= \frac{1}{2} f_1$, & $3$ & NA \\
& $\;\;e_2 f_2= \frac{1}{2} f_2$, $\;\;f_1f_2= e_1+ \gamma e_2$, $\; \gamma\in \F^*$. &&\\
$(2,2)_8$: & $\;\;e^2_1=e_1$, $\;\;e^2_2=e_2$, $\;\;e_1 f_1= \frac{1}{2} f_1$. & 2 & NA  \\
$(2,2)_9$: & $\;\;e^2_1=e_1$, $\;\;e^2_2=e_2$, $\;\;e_1 f_2= \frac{1}{2} f_2$, $\;\;e_2 f_1= \frac{1}{2} f_1$. & $2$ & NA  \\
$(2,2)_{10}$: & $\;\;e^2_1=e_1$, $\;\;e^2_2=e_2$, $\;\;e_1 f_1= f_1$, $\;\;e_2 f_2= \frac{1}{2} f_2$. & $2$& NA  \\
$(2,2)_{11}$: & $\;\;e^2_1=e_1$, $\;\;e^2_2=e_2$, $\;\;e_1 f_1= \frac{1}{2} f_1$, $\;\;e_2 f_1= \frac{1}{2} f_1$, 
$\;\;e_2 f_2= f_2$. & $2$ & NA  \\
$(2,2)_{12}$: & $\;\;e^2_1=e_1$, $\;\;e^2_2=e_2$, $\;\;e_1 f_1= f_1$.  & $2$ & A  \\
$(2,2)_{13}$: & $\;\;e^2_1=e_1$, $\;\;e^2_2=e_2$, $\;\;e_1 f_1= f_1$, $\;\;e_1 f_2= \frac{1}{2} f_2$. & $2$ & NA  \\
$(2,2)_{14}$: & $\;\;e^2_1=e_1$, $\;\;e^2_2=e_2$, $\;\;e_1 f_1= f_1$, $\;\;e_2 f_2= f_2$. & $2$ & A  \\
$(2,2)_{15}$: & $\;\;e^2_1=e_1$, $\;\;e^2_2=e_2$, $\;\;e_1 f_1= \frac{1}{2} f_1$, $\;\;e_2 f_1= \frac{1}{2}f_1$. & $2$ & NA \\
$(2,2)_{16}$: & $\;\;e^2_1=e_1$, $\;\;e^2_2=e_2$, $\;\;e_1 f_1= \frac{1}{2} f_1$,$\;\;e_1 f_2= \frac{1}{2} f_2$, $\;\;e_2 f_1= \frac{1}{2}f_1$. & $2$ & NA \\
$(2,2)_{17}$: & $\;\;e^2_1=e_1$. & $5$ & A  \\
$(2,2)_{18}$: & $\;\;e^2_1=e_1$, $\;\;f_1 f_2 =e_2$. & $4$ & A \\
$(2,2)_{19}$: & $\;\;e^2_1=e_1$, $\;\;e_2 f_2 =f_1$. & $3$ & A  \\
$(2,2)_{20}$: & $\;\;e^2_1=e_1$, $\;\;e_2 f_2 =f_1$, $\;\;f_1 f_2 =e_2$. & $2$& NA  \\
$(2,2)_{21}$: & $\;\;e^2_1=e_1$, $\;\;e_1 f_2 = \frac{1}{2} f_2$. & $3$ & NA  \\
$(2,2)_{22}$: & $\;\;e^2_1=e_1$, $\;\;e_1 f_2 = f_2$. & $3$ & A \\
$(2,2)_{23}$: & $\;\;e^2_1=e_1$, $\;\;e_1 f_1 = \frac{1}{2} f_1$, $\;\;e_1 f_2 = \frac{1}{2} f_2$. & $5$ & NA  \\
$(2,2)_{24}$: & $\;\;e^2_1=e_1$, $\;\;e_1 f_1 = \frac{1}{2} f_1$, $\;\;e_1 f_2 = \frac{1}{2} f_2$, $\;\;f_1 f_2 =e_1$. & $4$ & NA  \\
$(2,2)_{25}$: & $\;\;e^2_1=e_1$, $\;\;e_1 f_1 = \frac{1}{2} f_1$, $\;\;e_1 f_2 = \frac{1}{2} f_2$, $\;\;f_1 f_2 =e_2$. & $4$ & NA  \\
$(2,2)_{26}$: & $\;\;e^2_1=e_1$, $\;\;e_1 f_1 = \frac{1}{2} f_1$, $\;\;e_1 f_2 = \frac{1}{2} f_2$, 
$\;\;f_1 f_2 =e_1+e_2$. & $3$ & NA \\
$(2,2)_{27}$: & $\;\;e^2_1=e_1$, $\;\;e_1 f_1 = \frac{1}{2} f_1$, $\;\;e_1 f_2 = \frac{1}{2} f_2$, 
$\;\;e_2 f_2 =f_1$. & $3$& NA  \\
$(2,2)_{28}$: & $\;\;e^2_1=e_1$, $\;\;e_1 f_1 = \frac{1}{2} f_1$, $\;\;e_1 f_2 = \frac{1}{2} f_2$, $\;\;e_2 f_2 =f_1$, 
$\;\;f_1f_2=e_2$. & $2$ & NA \\
$(2,2)_{29}$: & $\;\;e^2_1=e_1$, $\;\;e_1 f_1 = \frac{1}{2} f_1$, $\;\;e_1 f_2 = f_2$. & $3$ & NA  \\
$(2,2)_{30}$: & $\;\;e^2_1=e_1$, $\;\;e_1 f_1 = f_1$, $\;\;e_1 f_2 = f_2$. & $5$ & A \\
$(2,2)_{31}$: & $\;\;e^2_1=e_1$, $\;\;e_1 f_1 = f_1$, $\;\;e_1 f_2 = f_2$, $\;\;f_1f_2=e_1$. & $4$ & NA  \\
$(2,2)_{32}$: & $\;\;e^2_1=e_1$, $\;\;e_1 e_2 = e_2$. & $5$& A  \\
$(2,2)_{33}$: & $\;\;e^2_1=e_1$, $\;\;e_1 e_2 = e_2$, $\;\;e_1 f_2 = \frac{1}{2} f_2$. & $3$& NA \\
$(2,2)_{34}$: & $\;\;e^2_1=e_1$, $\;\;e_1 e_2 = e_2$, $\;\;e_1 f_2 = f_2$. & $3$ & A  \\
$(2,2)_{35}$: & $\;\;e^2_1=e_1$, $\;\;e_1 e_2 = e_2$, $\;\;e_1 f_1 = \frac{1}{2} f_1$, $\;\;e_1 f_2 = \frac{1}{2} f_2$. & $5$ & NA  \\
$(2,2)_{36}$: &  $\;\;e^2_1=e_1$, $\;\;e_1 e_2 = e_2$, $\;\;e_1 f_1 = \frac{1}{2} f_1$, $\;\;e_1 f_2 = \frac{1}{2} f_2$, $\;\;f_1f_2=e_2$. & $4$ & NA \\
$(2,2)_{37}$: & $\;\;e^2_1=e_1$, $\;\;e_1 e_2 = e_2$, $\;\;e_1 f_1 = \frac{1}{2} f_1$, $\;\;e_1 f_2 = \frac{1}{2} f_2$, $\;\;e_2f_2=f_1$. & $3$& NA  \\
$(2,2)_{38}$: & $\;\;e^2_1=e_1$, $\;\;e_1 e_2 = e_2$, $\;\;e_1 f_1 = \frac{1}{2} f_1$, $\;\;e_1 f_2 = \frac{1}{2} f_2$, $\;\;e_2f_2=f_1$, 
& $2$ & NA  \\
&$\;\;f_1f_2=e_2$.&&\\
$(2,2)_{39}$: & $\;\;e^2_1=e_1$, $\;\;e_1 e_2 = e_2$, $\;\;e_1 f_1 = \frac{1}{2} f_1$, $\;\;e_1 f_2 = f_2$. & $3$ & NA  \\
$(2,2)_{40}$: & $\;\;e^2_1=e_1$, $\;\;e_1 e_2 = e_2$, $\;\;e_1 f_1 = f_1$, $\;\;e_1 f_2 = f_2$. & $5$ & A  \\
$(2,2)_{41}$: & $\;\;e^2_1=e_1$, $\;\;e_1 e_2 = e_2$, $\;\;e_1 f_1 = f_1$, $\;\;e_1 f_2 = f_2$, $\;\;f_1f_2=e_1$. & $4$ & NA  \\
$(2,2)_{42}$: & $\;\;e^2_1=e_1$, $\;\;e_1 e_2 = e_2$, $\;\;e_1 f_1 = f_1$, $\;\;e_1 f_2 = f_2$, 
$\;\;f_1f_2=e_2$. &$4$ & A  \\
$(2,2)_{43}$: & $\;\;e^2_1=e_1$, $\;\;e_1 e_2 = e_2$, $\;\;e_1 f_1 = f_1$, $\;\;e_1 f_2 = f_2$, $\;\;f_1f_2=e_1+e_2$. & $3$ & NA  \\
$(2,2)_{44}$: & $\;\;e^2_1=e_1$, $\;\;e_1 e_2 = e_2$, $\;\;e_1 f_1 = f_1$, $\;\;e_1 f_2 = f_2$, $\;\;e_2 f_2 = f_1$. & $3$ & A \\
$(2,2)_{45}$: & $\;\;e^2_1=e_1$, $\;\;e_1 e_2 = e_2$, $\;\;e_1 f_1 = f_1$, $\;\;e_1 f_2 = f_2$, $\;\;e_2 f_2 = f_1$, & $2$ & NA  \\
& $\;\;f_1f_2=e_2$.&&\\
$(2,2)^N_{46}$: & $\;\;e^2_1=e_2$. & $6$ & A  \\
$(2,2)^N_{47}$: & $\;\;e^2_1=e_2$, $\;\;f_1f_2=e_1$. & $4$ & NA  \\
$(2,2)^N_{48}$: & $\;\;e^2_1=e_2$, $\;\;f_1f_2=e_2$. & $5$ & A  \\
$(2,2)^N_{49}$: & $\;\;e^2_1=e_2$, $\;\;e_1f_2=f_1$. & $4$ & A  \\
$(2,2)^N_{50}$: & $\;\;e^2_1=e_2$, $\;\;e_1f_2=f_1$, $\;\;f_1f_2=e_2$. & $3$ & NA  \\
$(2,2)^N_{51}$: & $\;\;e^2_1=e_2$, $\;\;e_2f_2=f_1$. & $3$ & NA \\
$(2,2)_{52}$: & $\;\;e^2_1=e_1$, $\;\;e_1e_2= \frac{1}{2} e_2$. & $6$ & NA  \\
$(2,2)_{53}$: & $\;\;e^2_1=e_1$, $\;\;e_1e_2= \frac{1}{2} e_2$, $\;\;e_1f_2= \frac{1}{2} f_2$. & $4$ & NA \\
$(2,2)_{54}$: & $\;\;e^2_1=e_1$, $\;\;e_1e_2= \frac{1}{2} e_2$, $\;\;e_1f_2= \frac{1}{2} f_2$, $\;\;f_1f_2=e_2$. & $3$ & NA  \\
$(2,2)_{55}$: & $\;\;e^2_1=e_1$, $\;\;e_1e_2= \frac{1}{2} e_2$, $\;\;e_1f_2= \frac{1}{2} f_2$, $\;\;e_2f_2=f_1$. & $3$ & NA \\
$(2,2)_{56}$: & $\;\;e^2_1=e_1$, $\;\;e_1e_2= \frac{1}{2} e_2$, $\;\;e_1f_2= \frac{1}{2} f_2$, $\;\;e_2f_2=f_1$,
$\;\;f_1f_2=e_2$. & $2$& NA \\
$(2,2)_{57}$: & $\;\;e^2_1=e_1$, $\;\;e_1e_2= \frac{1}{2} e_2$, $\;\;e_1f_2= f_2$. & $4$ & NA \\
$(2,2)_{58}$: & $\;\;e^2_1=e_1$, $\;\;e_1e_2= \frac{1}{2} e_2$, $\;\;e_1f_1= \frac{1}{2} f_1$, $\;\;e_2f_2=f_1$. & $3$ & NA \\
$(2,2)_{59}$: & $\;\;e^2_1=e_1$, $\;\;e_1e_2= \frac{1}{2} e_2$, $\;\;e_1f_1= \frac{1}{2} f_1$, $\;\;e_2f_2=f_1$, 
$\;\;f_1f_2=e_2$. & $2$& NA  \\
$(2,2)_{60}$: & $\;\;e^2_1=e_1$, $\;\;e_1e_2= \frac{1}{2} e_2$, $\;\;e_1f_1= \frac{1}{2} f_1$, 
$\;\;e_1f_2= \frac{1}{2} f_2$. & $6$ & NA  \\
$(2,2)_{61}$: & $\;\;e^2_1=e_1$, $\;\;e_1e_2= \frac{1}{2} e_2$, $\;\;e_1f_1= \frac{1}{2} f_1$, $\;\;e_1f_2= f_2$. & $4$ & NA  \\
$(2,2)_{62}$: & $\;\;e^2_1=e_1$, $\;\;e_1e_2= \frac{1}{2} e_2$, $\;\;e_1f_1= \frac{1}{2} f_1$, $\;\;e_1f_2= f_2$, 
$\;\;f_1f_2=e_2$. & $3$ & NA \\
$(2,2)_{63}$: & $\;\;e^2_1=e_1$, $\;\;e_1e_2= \frac{1}{2} e_2$, $\;\;e_1f_1= \frac{1}{2} f_1$, $\;\;e_1f_2= f_2$, 
$\;\;e_2f_2= f_1$. & $3$& NA \\
$(2,2)_{64}$: & $\;\;e^2_1=e_1$, $\;\;e_1e_2= \frac{1}{2} e_2$, $\;\;e_1f_1= \frac{1}{2} f_1$, $\;\;e_1f_2= f_2$, $\;\;e_2f_2= f_1$, & $2$ & NA  \\
& $\;\;f_1f_2=e_2$.&&\\
$(2,2)_{65}$: & $\;\;e^2_1=e_1$, $\;\;e_1e_2= \frac{1}{2} e_2$, $\;\;e_1f_1= f_1$, $\;\;e_1f_2= \frac{1}{2} f_2$, $\;\;e_2f_2= f_1$.
& $3$ & NA \\
$(2,2)_{66}$: & $\;\;e^2_1=e_1$, $\;\;e_1e_2= \frac{1}{2} e_2$, $\;\;e_1f_1= f_1$, $\;\;e_1f_2= \frac{1}{2} f_2$, $\;\;e_2f_2= f_1$, & $2$ & NA \\
& $\;\;f_1f_2=e_2$.&&\\
$(2,2)_{67}$: & $\;\;e^2_1=e_1$, $\;\;e_1e_2= \frac{1}{2} e_2$, $\;\;e_1f_1= f_1$, $\;\;e_1f_2= f_2$. & $6$ & NA  \\
$(2,2)^N_{68}$: & $\;\;e_2f_2= f_1$. & $5$ & A  \\
$(2,2)^N_{69}$: & $\;\;e_2f_2= f_1$, $\;\;f_1f_2= e_1$. & $4$ & NA  \\
$(2,2)_{70}$: & $\;\;e_2f_2= f_1$, $\;\;f_1f_2= e_2$. & $3$ & NA  \\
$(2,2)^N_{71}$: & $\;\;f_1f_2= e_1$. & $6$ & A \\
$(2,2)^N_{72}$: & $\;\;e_ie_j=0$, $\;\;e_if_j= 0$, $\;\;f_if_j= 0$, $\;\; i, j \in \{1,2\}$.  & $8$ & A  \\
\hline
\end{longtable}


\begin{lemma} Every superalgebra $(2,2)_i\in \mathcal{JS}^{(2,2)}$ belongs to $\overline{\bigcup_{\gamma \in \mathbb{F}^*} D_{\gamma}^G}$ or $\overline{(2,2)_j^G}$ for $j\in \{1,3,5,8-16, 20, 28, 38, 45, 52, 56, 57, 59, 60,64, 66, 67\}$.
\end{lemma}
\begin{proof}
We will show that $(2,2)_{43}\in\overline{\bigcup\limits_{\gamma\in \mathbb{F}^{*}}D_t^{G}}$. Consider the ``variable'' change of basis $g(t)$ of $D_{\gamma}$, where the parameter $\gamma$ depends on $t$, namely: $\gamma(t)=1+t$ and $E_1=e_1+e_2$,
$E_2=te_2$, $F_1=f_1$ and 
$F_2=f_2$. So the curve $g(t)$ lies transversely 
to the orbits of $D_{\gamma}$, meaning that, for any $t\neq 0$, $g(t)\subset\bigcup\limits_{\gamma\in \mathbb{F}^{*}}D_{\gamma}^{G}$ and cuts $(2,2)_{43}^{G}$ in $t=0$. Thus $(2,2)_{43}\in\overline{\bigcup\limits_{\gamma\in \mathbb{F}^{*}}D_t^{G}}$. We will abuse the notation and will denote this fact by $D_t \to(2,2)_{43}$. Observe that in the same way we obtain $D_t \to(2,2)_{7}$ and $D_t \to(2,2)_{26}$ with the ``variable'' change of basis given in Table \ref{Table:deformacoes D_gamma}. 

\begin{table}[ht]
\begin{center}
\begin{spacing}{1.4}
\begin{tabular}{|c|llll|}
\hline
$D_\gamma \rightarrow (2,2)_{7}, \; \gamma(t) = t$ &	$E_1 = e_2$, & $E_2 =e_1$, & $F_1 = f_1$, & $F_2 = f_2$ \\
\hline
$D_\gamma \rightarrow (2,2)_{26}, \; \gamma(t) = t$ &	$E_1 = e_1$, & $E_2 =t e_2$, & $F_1 = f_1$, & $F_2 = f_2$ \\
\hline
$D_\gamma \rightarrow (2,2)_{43}, \; \gamma(t) = 1+t$ &	$E_1 = e_1+e_2$, & $E_2 =t e_2$, & $F_1 = f_1$, & $F_2 = f_2$ \\
\hline
\end{tabular}
\end{spacing}
\caption{\label{Table:deformacoes D_gamma}Change of basis of $D_{\gamma}$, where $\gamma$ depends on $t$}
\end{center}
\end{table}

On the other hand, Table \ref{table:change_bases_(2,2)} gives all possible essential deformations 
between Jordan superalgebras of type $(2,2)$ and the other deformations can be obtained by transitivity.

\begin{longtable}[H]{|l|llll|}
\caption{\label{table:change_bases_(2,2)}Deformations between Jordan Superalgebras of type $(2,2)$}\\
\hline
\multicolumn{1}{|c|}{$\mathcal{J} \rightarrow \mathcal{J}^\prime $ } & Change of basis & & & \\
\hline
\endhead
$(2,2)_{1} \rightarrow (2,2)_{17}$ & $E_1 = e_1$, & $E_2 = t e_2$, & $F_1 = f_1$, & $F_2 = f_2$ \\
\hline 
$(2,2)_{1} \rightarrow (2,2)_{32}$ & $E_1 = e_1 + e_2$, & $E_2 = t e_1$, & $F_1 = f_1$, & $F_2 = f_2$ \\
\hline 
$(2,2)_{2} \rightarrow (2,2)_{17}$ & $E_1 = e_2 $, & $E_2 = t e_1$, & $F_1 = f_1$, & $F_2 = f_2$ \\
\hline
$(2,2)_{2} \rightarrow (2,2)_{23}$ & $E_1 = e_1$, & $E_2 = t e_2$, & $F_1 = f_1$, & $F_2 = f_2$ \\
\hline
$(2,2)_{2} \rightarrow (2,2)_{35}$ & $E_1 = e_1 + e_2$, & $E_2 = t e_2$, & $F_1 = f_1$, & $F_2 = f_2$ \\
\hline
$(2,2)_{3} \rightarrow (2,2)_{2}$ & $E_1 = e_1$, & $E_2 = e_2$, & $F_1 = f_1$, & $F_2 = t f_2$ \\
\hline
$(2,2)_{3} \rightarrow (2,2)_{18}$ & $E_1 = e_2$, & $E_2 = t e_1$, & $F_1 = f_1$, & $F_2 = t f_2$ \\
\hline
$(2,2)_{3} \rightarrow (2,2)_{24}$ & $E_1 = e_1$, & $E_2 = te_2$, & $F_1 = f_1$, & $F_2 = f_2$ \\
\hline
$(2,2)_{3} \rightarrow (2,2)_{36}$ & $E_1 = e_1 + e_2$, & $E_2 = te_1$, & $F_1 = tf_1$, & $F_2 = f_2$ \\
\hline
$(2,2)_{4} \rightarrow (2,2)_{17}$ & $E_1 = e_2$, & $E_2 = t e_1$, & $F_1 = f_1$, & $F_2 = f_2$ \\
\hline
$(2,2)_{4} \rightarrow (2,2)_{30}$ & $E_1 = e_1$, & $E_2 = t e_2$, & $F_1 = f_1$, & $F_2 = f_2$ \\
\hline
$(2,2)_{4} \rightarrow (2,2)_{40}$ & $E_1 = e_1 + e_2$, & $E_2 = t e_2$, & $F_1 = f_1$, & $F_2 = f_2$ \\
\hline
$(2,2)_{5} \rightarrow (2,2)_{4}$ & $E_1 = e_1$, & $E_2 = e_2$, & $F_1 = f_1$, & $F_2 = tf_2$ \\
\hline
$(2,2)_{5} \rightarrow (2,2)_{18}$ & $E_1 = e_2$, & $E_2 = te_1$, & $F_1 = f_1$, & $F_2 = t f_2$ \\
\hline
$(2,2)_{5} \rightarrow (2,2)_{31}$ & $E_1 = e_1$, & $E_2 = te_2$, & $F_1 = f_1$, & $F_2 = f_2$ \\
\hline
$(2,2)_{5} \rightarrow (2,2)_{41}$ & $E_1 = e_1 + e_2$, & $E_2 = te_2$, & $F_1 = f_1$, & $F_2 = f_2$ \\
\hline
$(2,2)_{5} \rightarrow (2,2)_{42}$ & $E_1 = e_1 + e_2$, & $E_2 = t(e_1 - e_2)$, & $F_1 = 2t f_1$, & $F_2 = f_2$ \\
\hline
$(2,2)_{6} \rightarrow (2,2)_{23}$ & $E_1 = e_1$, & $E_2 = te_2$, & $F_1 = f_1$, & $F_2 = f_2$ \\
\hline
$(2,2)_{6} \rightarrow (2,2)_{40}$ & $E_1 = e_1 + e_2$, & $E_2 = t e_2$, & $F_1 = f_1$, & $F_2 = f_2$ \\
\hline
$(2,2)_{7} \rightarrow (2,2)_{6}$ & $E_1 = e_1$, & $E_2 = e_2$, & $F_1 = t f_1$, & $F_2 = f_2$ \\
\hline
$(2,2)_{7} \rightarrow (2,2)_{24}$ & $E_1 = e_2$, & $E_2 = t e_1$, & $F_1 = f_1$, & $F_2 = f_2$ \\
\hline
$(2,2)_{7} \rightarrow (2,2)_{25}$ & $E_1=e_1$, & $E_2=te_2$, & $F_1=t f_1$, & $F_2=f_2$\\
\hline
$(2,2)_{7} \rightarrow (2,2)_{42}$ & $E_1 = e_1 + e_2$, & $E_2 = t e_2$, & $F_1 = t f_1$, & $F_2 = f_2$ \\
\hline
$(2,2)_{8}\rightarrow(2,2)_{19}$ & $E_{1}=e_{2}$, & $E_{2}=te_{1}$, & $F_{1}=\frac{1}{2} tf_{1}$, & $F_{2}=f_{1}+f_{2}$\tabularnewline
\hline
$(2,2)_{8}\rightarrow(2,2)_{21}$ & $E_{1}=e_{1}$, & $E_{2}=te_{2}$, & $F_{1}=f_{2}$, & $F_{2}=f_{1}$.\tabularnewline
\hline 
 $(2,2)_{8}\rightarrow(2,2)_{33}$ & $E_{1}=e_{1}+e_{2}$, & $E_{2}=te_{2}$, & $F_{1}=f_{2}$, & $F_{2}=f_{1}$\tabularnewline
\hline 
$(2,2)_{9}\rightarrow(2,2)_{21}$ & $E_{1}=e_{1}$, & $E_{2}=te_{2}$, & $F_{1}=f_{1}$, & $F_{2}=f_{2}$\tabularnewline
\hline 
$(2,2)_{9} \rightarrow (2,2)_{37}$ & $E_1 = e_1 + e_2$, & $E_2 = t e_1 $, & $F_1 = \frac{1}{2} t f_2$, & $F_2 = f_1 + f_2$\\
\hline
$(2,2)_{10}\rightarrow(2,2)_{21}$ & $E_{1}=e_{2}$, & $E_{2}=te_{1}$, & $F_{1}=f_{1}$, & $F_{2}=f_{2}$\tabularnewline
\hline 
$(2,2)_{10}\rightarrow(2,2)_{22}$ & $E_{1}=e_{1}$, & $E_{2}=te_{2}$, & $F_{1}=f_{2}$, & $F_{2}=f_{1}$\tabularnewline
\hline
$(2,2)_{10} \rightarrow (2,2)_{39}$ & $E_1 = e_1 + e_2$, & $E_2 = t e_2$, & $F_1 = f_2$, & $F_2 = f_1$. \\
\hline
$(2,2)_{10} \rightarrow (2,2)^N_{51}$ & $E_1 = t e_1 + 2 t e_2$, & $E_2 = t^2 e_1 + 4 t^2 e_2$, & $F_1 = - t^2 f_1$, & $F_2 = f_1 + f_2$ \\
\hline

$(2,2)_{11} \rightarrow (2,2)_{21}$ &		$E_1 = e_1$, & $E_2 = t e_2$, & $F_1 = f_2$, & $F_2 = f_1$ \\		
\hline
$(2,2)_{11} \rightarrow (2,2)_{29}$	&	$E_1 = e_2$, & $E_2 = t e_1$, & $F_1 = f_1$, & $F_2 = f_2$ \\
\hline
$(2,2)_{11} \rightarrow (2,2)_{44}$ &		$E_1 = e_1 + e_2$, & $E_2 = t e_2$, & $F_1 =\frac{1}{2}t f_1 + tf_2$, & $F_2 = f_1+f_2$ \\
\hline
$(2,2)_{12}\rightarrow(2,2)_{19}$ & $E_{1}=e_{2}$, & $E_{2}=te_{1}$, & $F_{1}=tf_{1}$, & $F_{2}=f_{1}+f_{2}$\tabularnewline
\hline 
$(2,2)_{12}\rightarrow(2,2)_{22}$ & $E_{1}=e_{1}$, & $E_{2}=te_{2}$, & $F_{1}=f_{2}$, & $F_{2}=f_{1}$\tabularnewline
\hline 
$(2,2)_{12}\rightarrow(2,2)_{34}$ & $E_{1}=e_{1}+e_{2}$, & $E_{2}=te_{2}$, & $F_{1}=f_{2}$, & $F_{2}=f_{1}$\tabularnewline
\hline 
$(2,2)_{13}\rightarrow(2,2)_{19}$ & $E_{1}=e_{2}$, & $E_{2}=te_{1}$, & $F_{1}= t f_{1}+\frac{1}{2}t f_{2}$, & $F_{2}=f_{1}+f_{2}$\tabularnewline
\hline 
$(2,2)_{13}\rightarrow(2,2)_{29}$ & $E_{1}=e_{1}$, & $E_{2}=te_{2}$, & $F_{1}=f_{2}$, & $F_{2}=f_{1}$\tabularnewline
\hline 
$(2,2)_{13} \rightarrow (2,2)_{39}$ & $E_1 = e_1 + e_2$, & $E_2 = t e_2$, & $F_1 = f_2$, & $F_2 = f_1$ \\
\hline 
$(2,2)_{14}\rightarrow(2,2)_{22}$ & $E_{1}=e_{1}$, & $E_{2}=te_{2}$, & $F_{1}=f_{2}$, & $F_{2}=f_{1}$\tabularnewline
\hline 
$(2,2)_{14} \rightarrow (2,2)_{44}$ 	&	$E_1 = e_1 + e_2$, & $E_2 = t e_2$, & $F_1 = t f_2$, & $F_2 = f_1+ f_2$ \\
\hline
$(2,2)_{15}\rightarrow(2,2)_{21}$ & $E_{1}=e_{1}$, & $E_{2}=te_{2}$, & $F_{1}=f_{2}$, & $F_{2}=f_{1}$\tabularnewline
\hline
$(2,2)_{15}\rightarrow(2,2)_{34}$ & $E_{1}=e_{1}+e_{2}$, & $E_{2}=te_{2}$, & $F_{1}=f_{2}$, & $F_{2}=f_{1}$\tabularnewline
\hline 
$(2,2)_{15} \rightarrow (2,2)^N_{51}$ & $E_1 = -t e_1 + t e_2$, & $E_2 = t^2 e_1 + t^2 e_2$, & $F_1 = t^2 f_1$, & $F_2 = f_1 + f_2$\\
\hline
$(2,2)_{16} \rightarrow (2,2)_{21}$	&	$E_1 = e_2$, & $E_2 = t e_1$, & $F_1 = f_2$, & $F_2 = f_1$ \\
\hline		
$(2,2)_{16} \rightarrow (2,2)_{27}$ &	$E_1 = e_1$, & $E_2 = t e_2$, & $F_1 = \frac{1}{2} t f_1$, & $F_2 = f_1+ f_2$ \\
\hline
$(2,2)_{16} \rightarrow (2,2)_{39}$ &	$E_1 = e_1 + e_2$, & $E_2 = t e_2$, & $F_1 = f_2$, & $F_2 = f_1$ \\
\hline
$(2,2)_{17} \rightarrow (2,2)^N_{46}$ & $E_1 = te_1+e_2$, & $E_2 =t^2 e_1$, & $F_1 = f_1$, & $F_2 = f_2$ \\
\hline
$(2,2)_{18} \rightarrow (2,2)_{17}$ & $E_1 = e_1$, & $E_2 = t e_2$, & $F_1 = tf_1$, & $F_2 = t f_2$ \\
\hline
$(2,2)_{18} \rightarrow (2,2)^N_{48}$ & $E_1 = t e_1 + e_2$, & $E_2 = t^2 e_1$, & $F_1 = f_1$, & $F_2 = -t f_2$ \\
\hline
$(2,2)_{19} \rightarrow (2,2)_{17}$ & $E_1 = e_1$, & $E_2 = t e_2$, & $F_1 = t f_1$, & $F_2 = t f_2$
\\
\hline
$(2,2)_{19} \rightarrow (2,2)^N_{49}$ & $E_1 = te_1 + e_2$, & $E_2 = t^2 e_1$, & $F_1 = f_1$, & $F_2 = f_2$ \\
\hline
$(2,2)_{20}\rightarrow(2,2)_{18}$ & $E_{1}=e_{1}$, & $E_{2}=te_{2}$, & $F_{1}=f_{1}$, & $F_{2}=tf_{2}$\\
\hline 
$(2,2)_{20}\rightarrow(2,2)_{19}$ & $E_{1}=e_{1}$, & $E_{2}=e_{2}$, & $F_{1}=tf_{1}$, & $F_{2}=tf_{2}$\\
\hline 
$(2,2)_{20} \rightarrow (2,2)^N_{50}$ &	$E_1 =- t^2 e_1 + e_2$, & $E_2 = t^4 e_1$, & $F_1 = t f_1$, & $F_2 = t f_2$ \\
\hline
$(2,2)_{20}\rightarrow(2,2)_{70}$ & $E_{1}=t e_1$, & $E_{2}= e_2$, & $F_{1}= f_1$, & $F_{2}= f_2$\\
\hline 
$(2,2)_{21} \rightarrow (2,2)^N_{49}$ & $E_1 = t e_1 + e_2$, & $E_2 = t^2 e_1$, & $F_1 = f_2$, & $F_2 = f_1 + 2 \frac{1}{t} f_2$ \\
\hline
$(2,2)_{22} \rightarrow (2,2)^N_{49}$ & $E_1 = te_1 + e_2$, & $E_2 = t^2 e_1$, & $F_1 = tf_2$, & $F_2 = f_1 + f_2$ \\
\hline
$(2,2)_{23} \rightarrow (2,2)^N_{46}$ & $E_1 = t e_1 + e_2$, & $E_2 = t^2 e_1$, & $F_1 = f_1$, & $F_2 = f_2$ \\ 
\hline
$(2,2)_{24} \rightarrow (2,2)_{23}$ & $E_1 = e_1 $, & $E_2 = e_2$, & $F_1 = t f_1$, & $F_2 = f_2$ \\
\hline
$(2,2)_{24} \rightarrow (2,2)^N_{48}$ & $E_1 = t e_1 + e_2$, & $E_2 = t^2 e_1$, & $F_1 = t f_1$, & $F_2 = t f_2$ \\
\hline
$(2,2)_{25} \rightarrow (2,2)_{23}$ & $E_1 = e_1$, & $E_2 = t e_2$, & $F_1 = t f_1$, & $F_2 = t f_2$ \\
\hline
$(2,2)_{25} \rightarrow (2,2)^N_{48}$ & $E_1 = -t^2 e_1+e_2$, & $E_2 = t^4 e_1$, & $F_1 = t f_1$, & $F_2 = t f_2$ \\
\hline
$(2,2)_{26} \rightarrow (2,2)_{24}$ & $E_1 = e_1$, & $E_2 = \frac{1}{t} e_2$, & $F_1 = f_1$, & $F_2 = f_2$ \\
\hline
$(2,2)_{26} \rightarrow (2,2)_{25}$ & $E_1 = e_1$, & $E_2 = t e_2$, & $F_1 = t f_1$, & $F_2 = f_2$ \\
\hline
$(2,2)_{26} \rightarrow (2,2)^N_{47}$ & $E_1 = t (e_1 + e_2)$, & $E_2 = t^2 e_1$, & $F_1 = t f_1$, & $F_2 = f_2$ \\
\hline
$(2,2)_{27} \rightarrow (2,2)_{23}$ & $E_1 = e_1$, & $E_2 = t e_2$, & $F_1 = f_1$, & $F_2 = f_2$ \\
\hline
$(2,2)_{27} \rightarrow (2,2)^N_{49}$ & $E_1 = te_1 + e_2$, & $E_2 = t^2 e_1$, & $F_1 = -tf_2$, & $F_2 = \frac{1}{t} f_1 + f_2$ \\
\hline
$(2,2)_{28} \rightarrow (2,2)_{25}$ & 	$E_1 = e_1$, & $E_2 = t e_2$, & $F_1 = f_1$, & $F_2 = t f_2$ \\
\hline
$(2,2)_{28} \rightarrow (2,2)_{27}$ &		$E_1 = e_1$, & $E_2 = e_2$, & $F_1 = t f_1$, & $F_2 = t f_2$ \\
\hline
$(2,2)_{28} \rightarrow (2,2)^N_{50}$ &		$E_1 = -t^2 e_1 + e_2$, & $E_2 = t^4 e_1$, & $F_1 = - t f_1$, & $F_2 = - t f_2$ \\
\hline
$(2,2)_{28} \rightarrow (2,2)_{70}$	&	$E_1 =t e_1$, & $E_2 = e_2$, & $F_1 = f_1$, & $F_2 = f_2$ \\
\hline
$(2,2)_{29} \rightarrow (2,2)^N_{49}$ & $E_1 = te_1 + e_2$, & $E_2 = t^2 e_1$, & $F_1 = \frac{1}{2} t f_1+ t^2 f_2$, & $F_2 = f_1 + t f_2$ \\
\hline
$(2,2)_{30} \rightarrow (2,2)^N_{46}$ & $E_1 = t e_1 + t e_2$, & $E_2 = t^2 e_1$, & $F_1 = f_1$, & $F_2 = f_2$ \\
\hline
$(2,2)_{31} \rightarrow (2,2)_{30}$ & $E_1 = e_1$, & $E_2 = e_2$, & $F_1 = t f_1$, & $F_2 = f_2$ \\
\hline
$(2,2)_{31} \rightarrow (2,2)^N_{48}$ & $E_1 = t e_1 + e_2$, & $E_2 = t^2 e_1$, & $F_1 = t f_1$, & $F_2 = t f_2$ \\
\hline
$(2,2)_{32} \rightarrow (2,2)^N_{46}$ & $E_1 = t e_1 + e_2$, & $E_2 = - t^2 e_1$, & $F_1 = f_1$, & $F_2 = f_2$ \\
\hline
$(2,2)_{33} \rightarrow (2,2)^N_{49}$ & $E_1 = te_1 + e_2$, & $E_2 = - t^2 e_1$, & $F_1 = f_2$, & $F_2 = f_1 + 2 \frac{1}{t} f_2$ \\
\hline
$(2,2)_{34} \rightarrow (2,2)^N_{49}$ & $E_1 = t e_1 + e_2$, & $E_2 = - t^2 e_1$, & $F_1 = f_2$, & $F_2 = f_1 + \frac{1}{t} f_2$ \\
\hline
$(2,2)_{35} \rightarrow (2,2)^N_{46}$ & $E_1 = t e_1 + e_2$, & $E_2 =- t^2 e_1$, & $F_1 = f_1$, & $F_2 = f_2$ \\
\hline
$(2,2)_{36} \rightarrow (2,2)_{35}$ & $E_1 = e_1$, & $E_2 = e_2$, & $F_1 = tf_1$, & $F_2 = f_2$ \\
\hline
$(2,2)_{36} \rightarrow (2,2)^N_{48}$ & $E_1 = te_1+e_2$, & $E_2 = t^2 e_1 + 2t e_2$, & $F_1 =t f_1$, & $F_2 =f_2$ \\
\hline
$(2,2)_{37} \rightarrow (2,2)_{35}$ & $E_1 = e_1$, & $E_2 = t e_2$, & $F_1 = f_1$, & $F_2 = f_2$ \\
\hline
$(2,2)_{37} \rightarrow (2,2)^N_{49}$ & $E_1 = t(e_1 + e_2)$, & $E_2 = -t^2 e_1$, & $F_1 = t f_1+ \frac{1}{2}tf_2$, & $F_2 = f_2$ \\
\hline
$(2,2)_{38} \rightarrow (2,2)_{36}$	& $E_1 = e_1$, & $E_2 = t e_2$, & $F_1 = f_1$, & $F_2 =t f_2$ \\
\hline
$(2,2)_{38} \rightarrow (2,2)_{37}$	&	$E_1 = e_1$, & $E_2 = e_2$, & $F_1 = t f_1$, & $F_2 = t f_2$ \\
\hline
$(2,2)_{38} \rightarrow (2,2)^N_{50}$ &	$E_1 = t^2 e_1 + e_2$, & $E_2 = t^4 e_1 + 2t^2 e_2$, & $F_1 = t f_1$, & $F_2 = t  f_2$ \\
\hline
$(2,2)_{38} \rightarrow (2,2)_{70}$	&	$E_1 = te_1$, & $E_2 = e_2$, & $F_1 = f_1$, & $F_2 =f_2$ \\
\hline
$(2,2)_{39} \rightarrow (2,2)^N_{49}$ & $E_1 = te_1 + e_2$, & $E_2 = -t^2 e_1$, & $F_1 = f_1$, & $F_2 = - 2 \frac{1}{t} f_1 + f_2$ \\
\hline
$(2,2)_{40} \rightarrow (2,2)^N_{46}$ & $E_1 = t e_1 + e_2$, & $E_2 = -t^2 e_1$, & $F_1 = f_1$, & $F_2 = f_2$ \\
\hline
$(2,2)_{41} \rightarrow (2,2)_{40}$ & $E_1 = e_1$, & $E_2 = e_2$, & $F_1 = t f_1$, & $F_2 = f_2$ \\
\hline
$(2,2)_{41} \rightarrow (2,2)^N_{48}$ & $E_1 = t e_1 + e_2$, & $E_2 = t^2 e_1 + 2t e_2$, & $F_1 = f_2$, & $F_2 = t^2 f_1$ \\
\hline
$(2,2)_{42} \rightarrow (2,2)_{40}$ & $E_1 = e_1$, & $E_2 = e_2$, & $F_1 = t f_1$, & $F_2 = f_2$ \\
\hline
$(2,2)_{42} \rightarrow (2,2)^N_{48}$ & $E_1 = t e_1 + e_2$, & $E_2 = t^2 e_1 + 2t e_2$, & $F_1 = f_1$, & $F_2 = tf_2$ \\
\hline
$(2,2)_{43} \rightarrow (2,2)_{41}$ & $E_1 = e_1$, & $E_2 = \frac{1}{t} e_2$, & $F_1 = f_1$, & $F_2 = f_2$ \\
\hline
$(2,2)_{43} \rightarrow (2,2)_{42}$ & $E_1 = e_1$, & $E_2 = te_2$, & $F_1 = tf_1$, & $F_2 = f_2$ \\
\hline
$(2,2)_{43} \rightarrow (2,2)^N_{47}$ & $E_1 = t(e_1 + e_2)$, & $E_2 = -t^2 e_1$, & $F_1 = tf_1$, & $F_2 = f_2$ \\
\hline
$(2,2)_{44} \rightarrow (2,2)_{40}$ & $E_1 = e_1$, & $E_2 = te_2$, & $F_1 = f_1$, & $F_2 = f_2$ \\
\hline
$(2,2)_{44} \rightarrow (2,2)^N_{49}$ & $E_1 = t(e_1 + e_2)$, & $E_2 = - t^2 e_1$, & $F_1 = f_1$, & $F_2 = \frac{1}{t} f_2$ \\
\hline
$(2,2)_{45} \rightarrow (2,2)_{42}$ &	$E_1 = e_1$, & $E_2 = t e_2$, & $F_1 = f_1$, & $F_2 =t f_2$ \\
\hline
$(2,2)_{45} \rightarrow (2,2)_{44}$ &	$E_1 = e_1$, & $E_2 = e_2$, & $F_1 = t f_1$, & $F_2 = t f_2$ \\
\hline
$(2,2)_{45} \rightarrow (2,2)^N_{50}$ &	$E_1 = t^2 e_1 + t e_2$, & $E_2 = t^4 e_1 + 2 t^2 e_2$, & $F_1 = t f_1$, & $F_2 = t f_2$ \\
\hline
$(2,2)_{45} \rightarrow (2,2)_{70}$	&	$E_1 = te_1$, & $E_2 = e_2$, & $F_1 = f_1$, & $F_2 = f_2$ \\
\hline
$(2,2)^N_{47} \rightarrow (2,2)^N_{48}$ & $E_1 = t e_1 + e_2$, & $E_2 =- t^3 e_1$, & $F_1 = tf_1$, & $F_2 = -t^2 f_2$ \\
\hline
$(2,2)^N_{48} \rightarrow (2,2)^N_{46}$ & $E_1 = e_1$, & $E_2 = e_2$, & $F_1 = t f_1$, & $F_2 = f_2$ \\
\hline
$(2,2)^N_{48} \rightarrow (2,2)^N_{71}$ & $E_1 = e_2$, & $E_2 = t e_1$, & $F_1 = f_1$, & $F_2 = f_2$ \\
\hline
$(2,2)^N_{49} \rightarrow (2,2)^N_{46}$ & $E_1 = e_1$, & $E_2 = e_2$, & $F_1 = f_1$, & $F_2 = tf_2$ \\
\hline
$(2,2)^N_{49} \rightarrow (2,2)^N_{68}$ & $E_1 = e_2$, & $E_2 = t e_1$, & $F_1 = t f_1$, & $F_2 = f_2$ \\
\hline
$(2,2)^N_{50} \rightarrow (2,2)^N_{48}$ & $E_1 = t e_1$, & $E_2 = t^2 e_2$, & $F_1 = f_1$, & $F_2 = t^2 f_2$ \\
\hline
$(2,2)^N_{50} \rightarrow (2,2)^N_{49}$ & $E_1 = e_1$, & $E_2 = e_2$, & $F_1 = tf_1$, & $F_2 = tf_2$ \\
\hline
$(2,2)^N_{50} \rightarrow (2,2)^N_{69}$ & $E_1 = t e_2$, & $E_2 = t e_1$, & $F_1 = t f_1$, & $F_2 = f_2$ \\
\hline
$(2,2)^N_{51} \rightarrow (2,2)^N_{49}$ & $E_1 = t(e_1 + e_2)$, & $E_2 = t^2 e_2$, & $F_1 = f_1$, & $F_2 = \frac{1}{t} f_2$ \\
\hline
$(2,2)_{53} \rightarrow (2,2)^N_{68}$ & $E_1 = e_2$, & $E_2 = t e_1$, & $F_1 = \frac{1}{2} t f_2$, & $F_2 = f_1 + f_2$ \\
\hline
$(2,2)_{54} \rightarrow (2,2)_{53}$ & $E_1 = e_1$, & $E_2 = e_2$, & $F_1 = tf_1$, & $F_2 = f_2$ \\
\hline
$(2,2)_{54} \rightarrow (2,2)^N_{69}$ & $E_1 = \frac{t}{2} e_2$, & $E_2 = t e_1$, & $F_1 = \frac{t}{2} f_2$, & $F_2 = f_2- f_1$ \\
\hline
$(2,2)_{55} \rightarrow (2,2)_{53}$ & $E_1 = e_1$, & $E_2 = t e_2$, & $F_1 = f_1$, & $F_2 = f_2$ \\
\hline
$(2,2)_{56} \rightarrow (2,2)_{54}$	&	$E_1 = e_1$, & $E_2 = t e_2$, & $F_1 = f_1$, & $F_2 = t f_2$ \\
\hline
$(2,2)_{56} \rightarrow (2,2)_{55}$ &		$E_1 = e_1$, & $E_2 = e_2$, & $F_1 = t f_1$, & $F_2 = t f_2$ \\
\hline
$(2,2)_{56} \rightarrow (2,2)_{70}$	&	$E_1 = t e_1$, & $E_2 = e_2$, & $F_1 = f_1$, & $F_2 = f_2$ \\
\hline
$(2,2)_{57} \rightarrow (2,2)^N_{68}$ & $E_1 = e_2$, & $E_2 = t e_1$, & $F_1 = t f_2$, & $F_2 = f_1 + f_2$ \\
\hline
$(2,2)_{58} \rightarrow (2,2)_{53}$ & $E_1 = e_1$, & $E_2 = t e_2$, & $F_1 = f_2$, & $F_2 = f_1$ \\
\hline
$(2,2)_{59} \rightarrow (2,2)_{54}$	&	$E_1 = e_1$, & $E_2 = t^2 e_2$, & $F_1 = t f_2$, & $F_2 = t f_1$ \\
\hline
$(2,2)_{59} \rightarrow (2,2)_{58}$	&	$E_1 = e_1$, & $E_2 = e_2$, & $F_1 = t f_1$, & $F_2 = t f_2$ \\
\hline
$(2,2)_{59} \rightarrow (2,2)_{70}$	&	$E_1 = te_1$, & $E_2 = e_2$, & $F_1 = f_1$, & $F_2 = f_2$ \\
\hline
$(2,2)_{61} \rightarrow (2,2)^N_{68}$ & $E_1 = e_2$, & $E_2 = t e_1$, & $F_1 = \frac{1}{2} t f_1 + t f_2$, & $F_2 = f_1 + f_2$ \\
\hline
$(2,2)_{62} \rightarrow (2,2)_{61}$ & $E_1 = e_1$, & $E_2 = e_2$, & $F_1 = f_1$, & $F_2 = tf_2$ \\
\hline
$(2,2)_{62} \rightarrow (2,2)^N_{69}$ & $E_1 = -2t e_2$, & $E_2 = t e_1$, & $F_1 = t f_2$, & $F_2 =2( f_1 + f_2)$ \\
\hline
$(2,2)_{63} \rightarrow (2,2)_{61}$ & $E_1 = e_1$, & $E_2 = t e_2$, & $F_1 = f_1$, & $F_2 = f_2$ \\
\hline
$(2,2)_{64} \rightarrow (2,2)_{62}$	&	$E_1 = e_1$, & $E_2 = t e_2$, & $F_1 = f_1$, & $F_2 = t f_2$ \\
\hline
$(2,2)_{64} \rightarrow (2,2)_{63}$	&	$E_1 = e_1$, & $E_2 = e_2$, & $F_1 = t f_1$, & $F_2 = t f_2$ \\
\hline
$(2,2)_{64} \rightarrow (2,2)_{70}$	&	$E_1 = t e_1$, & $E_2 = e_2$, & $F_1 = f_1$, & $F_2 = f_2$ \\
\hline
$(2,2)_{65} \rightarrow (2,2)_{61}$ & $E_1 = e_1$, & $E_2 = te_2$, & $F_1 = f_2$, & $F_2 = f_1$ \\
\hline
$(2,2)_{66} \rightarrow (2,2)_{62}$	&	$E_1 = e_1$, & $E_2 = t e_2$, & $F_1 = f_1$, & $F_2 = t f_2$ \\
\hline
$(2,2)_{66} \rightarrow (2,2)_{65}$	&	$E_1 = e_1$, & $E_2 = e_2$, & $F_1 = t f_1$, & $F_2 = t f_2$ \\	
\hline
$(2,2)_{66} \rightarrow (2,2)_{70}$	&	$E_1 = t e_1$, & $E_2 = e_2$, & $F_1 = f_1$, & $F_2 = f_2$ \\
\hline
$(2,2)^N_{69} \rightarrow (2,2)^N_{68}$ & $E_1 = e_1$, & $E_2 = e_2$, & $F_1 = t f_1$, & $F_2 = t f_2$ \\
\hline
$(2,2)^N_{69} \rightarrow (2,2)^N_{71}$ & $E_1 = e_1$, & $E_2 = t e_2$, & $F_1 = f_1$, & $F_2 = f_2$ \\
\hline
$(2,2)_{70} \rightarrow (2,2)^N_{69}$ &	$E_1 = t^2 e_2$, & $E_2 = e_1+ e_2$, & $F_1 = t f_1$, & $F_2 = t f_2$ \\
\hline
$(2,2)_{i} \rightarrow (2,2)^N_{72}$ &	$E_1 = t e_1$, & $E_2 = t e_2$, & $F_1 = t f_1$, & $F_2 = t f_2$ \\
\hline
\end{longtable}
\end{proof}

\begin{lemma}\label{Lema:familianaocontida} The infinite union $\bigcup_{\gamma \in \mathbb{F}^*} D_{\gamma}^G$ is not contained in the Zariski closure of the orbits of any superalgebra in $\mathcal{JS}^{(2,2)}$.
\end{lemma}

\begin{proof} Suppose, for the sake of contradiction, that $\bigcup_{\gamma \in \mathbb{F}^*} D_{\gamma}^G \subseteq \overline{(2,2)_j^G}$ for some $1 \leq j \leq 72$. Hence, $D_{\gamma}^G \subseteq \overline{(2,2)_j^G}$ for all $\gamma \in \mathbb{F}^*$, which means $(2,2)_j \to D_{\gamma}$.
Since $\dim \Aut(D_{\gamma}) = 3$, it follows from Lemma \ref{Lema:nondeformation}(\ref{dimAut}) that $j \in \{8-16, 20, 28, 38, 45, 56, 59, 64, 66\}$. Now, since $(D_{\gamma})_0 = \mathbb{F}e_1 \oplus \mathbb{F}e_2$, which is a semisimple Jordan algebra of dimension $2$ and thus rigid, it follows from Lemma \ref{Lema:nondeformation}(\ref{partepar}) that $j \in \{8, \dots, 16\}$. However, $((2,2)_j)_1^2 = 0$ for all $j \in \{8, \dots, 16\}$, while $\dim (D_{\gamma})_1^2 = 1$, which contradicts Lemma \ref{Lema:nondeformation}(\ref{dimJ_1^2}).
\end{proof}

\begin{corollary} $\overline{\bigcup_{\gamma \in \mathbb{F}^*} D_{\gamma}^G}$ determines a component of the variety $\mathcal{JS}^{(2,2)}$.
\end{corollary}
\begin{proof}
 It follows from Lemma \ref{Lema:familianaocontida} that it is sufficient to show that all superalgebras in $\overline{\bigcup_{\gamma \in \mathbb{F}^*} D_{\gamma}^G}$ belong to the same component. Indeed, $\overline{\bigcup_{\gamma \in \mathbb{F}^*} D_{\gamma}^G}$ is the product of two irreducible varieties, $D_{\gamma}^G$ and $\mathbb{F}$, and is therefore irreducible.
\end{proof}

\begin{lemma}Superalgebras $(2,2)_j$ for $j\in \{1,3,5,8-16, 20, 28, 38, 45, 52, 56, 57, 59, 60,64, 66, 67\}$ are rigid.
\end{lemma}
\begin{proof} 
Examples \ref{ex:(2,2)_3_rigid}, \ref{ex:(2,2)_5_rigid}, and \ref{ex:(2,2)_1_rigid} show that the superalgebras $(2,2)_3$, $(2,2)_5$ and $(2,2)_1$ are rigid. Analogous calculations show that $(2,2)_{52}$, $(2,2)_{60}$, and $(2,2)_{67}$ are also rigid. 

Consider the superalgebras $(2,2)_j$ where $j \in \{8, \dots, 16, 20, 28, 38, 45, 56, 59, 64, 66\}$. All of them have automorphism group dimension equal to $2$, and thus $\dim (2,2)_j^G = \dim \overline{\bigcup_{\gamma \in \mathbb{F}^*} D_{\gamma}^G}$. Suppose that $(2,2)_j$ belongs to the component $\overline{\bigcup_{\gamma \in \mathbb{F}^*} D_{\gamma}^G} = \mathcal{C}$. Then we must have $\overline{(2,2)_j^G} = \mathcal{C}$, meaning $(2,2)_j^G$ is a dense subset of $\mathcal{C}$, and $\bigcup_{\gamma \in \mathbb{F}^*} D_{\gamma}^G$ is open in $\mathcal{C}$. Therefore, there exists $J \in (2,2)_j^G$ such that $J \in \bigcup_{\gamma \in \mathbb{F}^*} D_{\gamma}^G$, which implies there exists $\gamma \in \mathbb{F}^*$ such that $(2,2)_j \simeq D_{\gamma}$. This is a contradiction. Thus, $(2,2)_j^G \not\subseteq \overline{\bigcup_{\gamma \in \mathbb{F}^*} D_{\gamma}^G}$. 

Let us also verify that $(2,2)_{57}^G$ is not contained in the component $\overline{\bigcup_{\gamma \in \mathbb{F}^*} D_{\gamma}^G}$. Suppose it is; then, by Lemma \ref{Lema:curva-Mazzola}, there exists a curve $g(t) \subset \bigcup_{\gamma \in \mathbb{F}^*} D_{\gamma}^G$ given by  
\[
E_1^t = a_{11}(t)e_1 + a_{21}(t)e_2, \quad 
E_2^t = a_{12}(t)e_1 + a_{22}(t)e_2, \\
\]
\[
F_1^t = a_{33}(t)f_1 + a_{43}(t)f_2, \quad 
F_2^t = a_{34}(t)f_1 + a_{44}(t)f_2, 
\]
such that for some $t_0 \in \mathbb{F}$, we have $g(t_0) = (2,2)_{57}$.

Observe that in $(2,2)_{57}$, we have $e_1f_1 = f_1$ and $e_1f_2 = 0$. For $g(t)$, we have 
\[
E_i^t F_j^t = \frac{1}{2}(a_{1i}(t) + a_{2i}(t))F_j^t, \quad \text{for } i,j = 1,2.
\]
This shows that for every \(t \in \mathbb{F}\), \(E_1^t F_1^t\) and \(E_1^t F_2^t\) are scalar multiples of \(F_1^t\) and \(F_2^t\), respectively, with the same scalar for both products. Thus, there does not exist \(t_0 \in \mathbb{F}\) such that $g(t_0) = (2,2)_{57}$.

On the other hand, as a consequence of the information given in Table \ref{table:JSA_(2,2)}, the lowest dimension of an automorphism group is $2$ then it follows from Lemma \ref{Lema:nondeformation}(\ref{dimAut}) that all deformations of the superalgebras $(2,2)_j$, for $j \in \{8, \dots, 16, 20, 28, 38, 45,56, 59, 64, 66\}$ are trivial. Hence, they are rigid superalgebras. 
Finally, consider the superalgebra $(2,2)_{57}$. 
According to Table \ref{table:non_deformations_(2,2)},
$(2,2)_j\not \to (2,2)_{57}$
for $j\in \{1,\dots, 51 \}\cup\{68,\dots, 72\}$ by Lemma \ref{Lema:nondeformation}(\ref{partepar}), also $(2,2)_j\not \to (2,2)_{57}$ for $j\in \{54,56,58,59,62,64,66\}$ by Lemma \ref{Lema:nondeformation}(\ref{esquecimento}), $(2,2)_j\not \to (2,2)_{57}$  for $j\in \{55,63,65\}$ by Remark \ref{Remark:como algebras} and the information given in \cite{tesejenny}. Finally, $(2,2)_{j}\not \to (2,2)_{57}$ for $j\in \{52,53,60,61,67\}$ by Lemma \ref{Lema:nondeformation}(\ref{dimAut}) proving that $(2,2)_{57}$ is a rigid Jordan superalgebra.
\end{proof}

To determine the associated Hasse diagram and describe the irreducible components of the variety $\mathcal{JS}^{(2,2)}$, we present Table \ref{table:non_deformations_(2,2)}, which confirms the non-existence of deformations $(2,2)_i \not\rightarrow (2,2)_j$. For brevity, we denote this by $i \not\rightarrow j$.

\begin{longtable}[H]{| l | c |}
\caption{\label{table:non_deformations_(2,2)} Non-deformations for Jordan Superalgebras of type $(2, 2)$} 
\endhead
\hline 
\multicolumn{1}{|c|}{$\mathcal{J} \not\rightarrow \mathcal{J}^\prime $ } & Reason \\
\hline
$ i \not\rightarrow j$, $i \in \{1, \cdots, 16\}$ and 
 $j \in \{52, \cdots, 67\}$; & \\ 
$ i \not\rightarrow j$, for $i \in \{17, \cdots, 31\}$ and 
 $j \in \{1, \cdots, 16, 32, \cdots, 45, 52, \cdots, 67\}$; & 
 \\ 
$ i \not\rightarrow j$, for $i \in \{32, \cdots, 45\}$ and 
$j \in \{1, \cdots, 31, 52, \cdots, 67\}$; &  $\mathcal{J}_0 \not\rightarrow \mathcal{J}_0^\prime$
\\ 
$ i \not\rightarrow j$, for $i \in \{46, \cdots, 51\}$ and 
 $j \in \{1, \cdots, 45, 52, \cdots, 67\}$; & 
\\ 
$ i \not\rightarrow j$, for $i \in \{52, \cdots, 67\}$ and 
 $j \in \{1, \cdots, 51\}$; & 
\\ 
 $ i \not\rightarrow j$, for $i \in \{68, \cdots, 72\}$ and 
 $j \in \{1, \cdots, 67\}$. & 
\\ 
\hline
\hline
 $ i \not\rightarrow j$, for $i \in \{8, \cdots,16\}$ and &
 \\
 $j \in \{3, 5, 7, 18, 24, 25, 26, 31, 36,  41, 42, 43, 47, 48, 50, 69, 70, 71\}$; &
 \\
 $ i \not\rightarrow j$, for $i \in \{ 19, 21, 22, 27, 29 \}$ and 
 $j \in \{18, 24, 25, 31, 47, 48, 69, 71 \}$; &
 \\ 
 $ i \not\rightarrow j$, for $i \in \{ 33, 34, 37, 39, 44\}$ and  $j \in \{ 36, 41, 42, 47, 48, 69, 71 \}$; & $\operatorname{ab}(\mathcal{J}) \not\rightarrow \operatorname{ab}(\mathcal{J}^\prime)$
 \\ 
  $ 51 \not\rightarrow j$, for 
 $j \in \{47, 48, 69, 71\}$; & \\
  $ i \not\rightarrow j$, for $i \in \{ 55, 58, 63, 65\}$ and 
 $j \in \{ 69, 71\}$; &
 \\ 
$ i \not\rightarrow j$, for $i \in \{1, 2, 4, 6, 49\}$ and 
 $j \in \{48, 71 \}$; &
 \\
 $ i \not\rightarrow 71$, for $i \in \{17, 23, 30, 32, 35, 40, 53, 57, 61, 68\}$. 
 &
 \\  
\hline
\hline
$2  \not\rightarrow \{30, 32, 40 \}$; 
$ \; \; 4 \not\rightarrow 32 $; 
$ 6 \not\rightarrow \{17, 30, 32, 40 \} $;& 
\\
$8 \not\rightarrow \{1, 2, 4, 6, 22, 23, 27, 30, 32, 34, 35, 37, 39, 40, 44\} $;& 
\\
$ 9 \not\rightarrow \{   1,2,4,6, 17,19, 22, 23, 27, 29, 30, 32, 33, 34, 39, 40, 44, 51 \}  $; &
\\
$ 10 \not\rightarrow \{ 1, 2, 4, 6, 17, 19, 23, 27, 29, 30, 33, 34, 35, 37, 40, 44  \} $; &
\\
$ 11 \not\rightarrow \{ 1, 2, 4, 6, 17, 19, 22, 23, 27, 30, 32, 34, 35, 37, 39 \} $; 
$ \; \; 12 \not\rightarrow \{ 1, 4, 40, 44 \} $;& 
\\
$ 13 \not\rightarrow \{1, 2, 4, 6, 21, 22, 23, 27, 30, 32, 34, 35, 37, 40, 44\} $; 
& $\mathcal{J} \not\rightarrow \mathcal{J}^\prime$\\
$ 14 \not\rightarrow \{ 1, 4, 17, 19, 30, 34 \} $;& 
as algebras\\
$15 \not\rightarrow \{ 1, 2, 4, 6, 17, 19, 22, 23, 27 ,29, 32, 33, 35, 37, 39, 40, 44 \} $;&
\\
$16 \not\rightarrow \{ 1, 2, 4, 6, 17, 19, 22, 29, 30, 32, 33, 34, 35, 37, 40, 44, 51\} $;& 
\\
$21, 29 \not\rightarrow \{ 17, 23, 30 \} $; 
$ \; \; 22, 27 \not\rightarrow \{17, 30 \} $; 
$\; \; 33, 39 \not\rightarrow \{ 32, 35, 40 \} $;&
\\ 
$34, 37 \not\rightarrow \{ 32, 40\} $;
$\; \; 44 \not\rightarrow 32 $; 
$\; \; 55 \not\rightarrow \{ 52, 57, 60, 61, 67 \}$; &
\\
$i \not\rightarrow \{ 52, 60, 67 \}$, for $i\in \{53, 57, 58, 61 \}$; 
$\; \; 63, 65 \not\rightarrow \{ 52, 53, 57, 60, 67 \} $.& 
\\
\hline
\hline
$ 3 \not\rightarrow \{1, 4, 6, 30, 31, 32, 40, 41, 42, 49, 68, 69 \} $;&
\\
$ 5 \not\rightarrow \{1, 2, 6, 23, 24, 25, 32, 35, 36, 49, 68, 69\} $; & 
\\
$ 7 \not\rightarrow \{ 1, 2, 4, 17, 18, 30, 31, 32, 35, 36, 49, 68, 69 \}$; & 
\\
$20 \not\rightarrow \{ 21, \cdots, 27, 29, 30, 31, 51 \} $;
$\; \; 24, 25  \not\rightarrow \{17, 30 \} $; & 
\\
$26 \not\rightarrow \{17, 18, 30, 31, 49, 68, 69 \} $; 
$\; \; 28 \not\rightarrow \{17, 18, 19, 21, 22, 29, 30, 31, 51  \} $; &
$\mathcal{F}(\mathcal{J}) \not\rightarrow \mathcal{F}(\mathcal{J}^\prime) $
\\
$31\not\rightarrow \{ 17, 23 \} $; $\; \; 36 \not\rightarrow \{32, 40 \}$;
$\; \; 38 \not\rightarrow \{32, 33, 34, 39, 40, 41, 42, 43, 44, 51 \} $; 
&
\\
$41 \not\rightarrow \{ 32, 35, 68\} $; $ \; \; 42 \not\rightarrow \{68\} $;
$\; \; 43  \not\rightarrow \{32, 35, 36, 49, 68, 69 \} $; 
&
\\
$45 \not\rightarrow \{ 33, \cdots, 37, 39, 51 \} $;
$\; \; 47 \not\rightarrow \{ 33 \} $;
$\; \; 54 \not\rightarrow \{ 52, 57, 60, 61, 67 \} $; 
&
\\
$56 \not\rightarrow \{ 52, 57, 58, 60, 61, 62, 63, 65, 67\} $; $ \;58 \not\rightarrow \{57 \}$; 
& 
\\
$59 \not\rightarrow \{ 52, 55, 57, 58, 60, 61, 62, 63, 65, 67 \} $;
$\; \; 62  \not\rightarrow \{ 52, 53, 57, 60, 67 \} $;
& 
\\
$ 64 \not\rightarrow \{ 52, \cdots, 55, 57, 58, 60, 65, 67 \} $;
$ \; \; 66 \not\rightarrow \{ 52, \cdots, 55, 57, 58, 60, 63, 67 \} $.&
\\
\hline
$50 \not\rightarrow 47$. &
$\dim ( \mathcal{J}^{2})_0 < \dim ( (\mathcal{J}^\prime)^{2})_0 $ 
\\
\hline
\hline
$ 1 \not\rightarrow \{ 30, 40, 68 \}$; \;\; $ 8 \not\rightarrow 29$;
$\;\; 18\not\rightarrow \{ 30, 68 \}$;& 
\\
$i \not\rightarrow 30$ for $i \in \{ 12, 19, 15 \}$ ; $ \;\; 47\not\rightarrow 68$ 
$\;\; 58 \not\rightarrow 61$.  &
$\dim ( \mathcal{J}^{2})_1 < \dim ( (\mathcal{J}^\prime)^{2})_1 $
\\
\hline
\hline
$ i  \not\rightarrow 68 $, for $i \in \{ 2, 4, 6, 24, 25, 31, 36 \}$; &
\\
$D_{\gamma}\not \to i$, for $i\in \{8,\cdots,16,19,\cdots,22,27,28,29,33,34,37,38,39,44,$ & ``general basis"\\
$45,49,50,51,53,\cdots,59,61,\cdots,66,68,69,70\}$ & \\
\hline
\end{longtable}

Notice that the superalgebras $(2,2)_i$ are not deformations of $(2,2)_{68}$, for $i \in \{ 2, 4, 6, 24, 25, 31,36 \}$ according to a criterion we call the ``\textit{general basis}''.
In fact, suppose that $(2,2)_i \to (2,2)_{68}$,
for $i\in \{2,4,6,24,25,31,36\}$. Then, for each $i$ there exists a
parameterized basis, namely:
\begin{align*}
E_1^t & = a_{11}(t)e_1 + a_{21}(t)e_2, & 
E_2^t &= a_{12}(t)e_1 + a_{22}(t)e_2, \\
F_1^t & = a_{33}(t) f_1+a_{43}(t) f_2, &
F_2^t & = a_{34}(t) f_1+a_{44}(t) f_2, \\
\end{align*}
such that for $t = 0$ we obtain $68$. Observe that
in $68$ we have $e_2f_2=f_1$ and 
\begin{align*}
E_2^t F_2^t & =\frac{1}{2} a_{12}(t)F_2^t, \,\text{ for 
}\,i=2,24,25,36.\\
E_2^t F_2^t &= a_{12}(t)F_2^t, \,\text{ for }\,i=4,31.\\
E_2^t F_2^t & = \frac{1}{2} (a_{12}(t)+a_{22}(t)) F_2^t, \,\text{ for 
}\,i=6.
\end{align*}
This shows that in all cases, \( E_2^tF_2^t \) is a multiple of \( F_2^t \), and by the linear independence of the basis, it follows that \( E_2^tF_2^t \neq F_1^t \) for all \( t\in \mathbb{F} \), which implies that \((2,2)_i \not \to (2,2)_{68}\).

According to the above information, we deduce 
the principal result of this work. 

\begin{theorem}\label{Thm: principal}
The variety $\mathcal{JS}^{(2,2)}$ has $25$ irreducible components, one of them is given by $\overline{ \bigcup_\gamma D^G_{\gamma}}$, the others are given by the Zariski closure of orbits of rigid superalgebras and they are the following:
\begin{longtable}[H]{rcl}

$\overline{ \bigcup_\gamma D^G_{\gamma}} $ & $=$ & 
$ \{ 
D_\gamma, {\color{red}(2,2)_{2}, (2,2)_{4}},    (2,2)_{6}, (2,2)_{7},{\color{red}(2,2)_{17}, (2,2)_{18}}, (2,2)_{23}, (2,2)_{24}, (2,2)_{25},  
$\\ 
&& $(2,2)_{26}, {\color{red}(2,2)_{30}, (2,2)_{31}, (2,2)_{32}, (2,2)_{35}, (2,2)_{36}}, 
(2,2)_{40}, (2,2)_{41}, (2,2)_{42},$ \\
&& $(2,2)_{43}, (2,2)^N_{46}, (2,2)^N_{47}, 
(2,2)^N_{48}, (2,2)^N_{71}, (2,2)^N_{72}
\} $.\\ 
$\overline{ (2,2)^G_{1} } $ & $=$ & 
$ \{ 
(2,2)_{1}, (2,2)_{17}, (2,2)_{32}, (2,2)^N_{46}, (2,2)^N_{72}
\} $. \\ 
$\overline{ (2,2)^G_{3} } $ & $=$ & 
$ \{ 
(2,2)_{2}, (2,2)_{3}, (2,2)_{17}, (2,2)_{18}, (2,2)_{23}, (2,2)_{24}, {\color{red}(2,2)_{25}}, (2,2)_{35}, (2,2)_{36}, 
$\\ 
&& $
(2,2)^N_{46}, {\color{red} (2,2)^N_{47} }, (2,2)^N_{48},  (2,2)^N_{71}, (2,2)^N_{72}
\} $. \\ 
$\overline{ (2,2)^G_{5} } $ & $=$ & 
$ \{ 
(2,2)_{4}, (2,2)_{5}, (2,2)_{17}, (2,2)_{18}, (2,2)_{30}, (2,2)_{31}, (2,2)_{40}, (2,2)_{41}, (2,2)_{42}$,\\ 
&& $(2,2)^N_{46}, {\color{red} (2,2)^N_{47} }, (2,2)^N_{48}, (2,2)^N_{71}, (2,2)^N_{72}
\} $. \\ 
$\overline{ (2,2)^G_{8} } $ & $=$ & 
$\{  
(2,2)_{8}, (2,2)_{17}, (2,2)_{19}, (2,2)_{21}, (2,2)_{33}, (2,2)^N_{46}, (2,2)^N_{49}, {\color{red}(2,2)^N_{51}}, 
$\\ 
&& $ (2,2)^N_{68}, (2,2)^N_{72}\}$. 
\\
$\overline{ (2,2)^G_{9} } $ & $=$ & 
$ \{ 
(2,2)_{9}, (2,2)_{21}, (2,2)_{35}, (2,2)_{37}, (2,2)^N_{46}, (2,2)^N_{49}, (2,2)^N_{68}, (2,2)^N_{72} 
\} $.\\
&&\\
$\overline{ (2,2)^G_{10} } $ & $=$ & 
$ \{ (2,2)_{10},  (2,2)_{ 21}, (2,2)_{22}, {\color{red} (2,2)_{32}}, (2,2)_{39}, (2,2)^N_{46}, (2,2)^N_{49}, (2,2)^N_{51}, (2,2)^N_{68},$\\ 
&& $ (2,2)^N_{72}  \} $.\\
 $\overline{ (2,2)^G_{11} } $ & $=$ & 
 $ \{ 
 (2,2)_{11}, (2,2)_{21}, (2,2)_{29}, {\color{red} (2,2)_{33}}, (2,2)_{40} , (2,2)_{44}, (2,2)^N_{46}, (2,2)^N_{49}, {\color{red} (2,2)_{ 51}},  
 $\\ 
&& $ (2,2)^N_{68}, (2,2)^N_{72} \} $.\\
 $\overline{ (2,2)^G_{12} } $ & $=$ & $ \{ 
(2,2)_{12},  (2,2)_{17}, (2,2)_{19}, (2,2)_{22}, {\color{red} (2,2)_{32}}, (2,2)_{34}, (2,2)^N_{46}, (2,2)^N_{49},$\\ 
&& $ (2,2)^N_{68}, (2,2)^N_{72}
\} $.\\
 $\overline{ (2,2)^G_{13} } $ & $=$ & 
 $ \{ (2,2)_{13}, (2,2)_{17}, (2,2)_{19}, (2,2)_{29}, {\color{red}(2,2)_{33}}, (2,2)_{39}, (2,2)^N_{46}, (2,2)^N_{49}, {\color{red} (2,2)_{ 51}},$\\ 
&& $(2,2)^N_{68}, (2,2)^N_{72}  \}  $. \\ 
 $\overline{ (2,2)^G_{14} } $ & $=$ &
$ \{ (2,2)_{14}, (2,2)_{22}, {\color{red} (2,2)_{32}},  (2,2)_{40}, (2,2)_{44}, (2,2)^N_{46}, (2,2)^N_{49}, (2,2)^N_{68}, (2,2)^N_{72}
\}  $. \\ 
&&\\
$\overline{ (2,2)^G_{15} } $ & $=$ & 
$ \{ (2,2)_{15},  (2,2)_{21}, (2,2)_{34}, (2,2)^N_{46}, (2,2)^N_{49}, (2,2)^N_{51}, (2,2)^N_{68}, (2,2)^N_{72}
\}  $. \\ 
&&\\
$\overline{ (2,2)^G_{16} } $ & $=$ & 
$ \{ (2,2)_{16}, (2,2)_{ 21}, (2,2)_{23}, (2,2)_{27}, (2,2)_{39}, (2,2)^N_{46}, (2,2)^N_{49}, (2,2)^N_{68}, (2,2)^N_{72} \}  $. \\ 
 &&\\
$\overline{ (2,2)^G_{20} }$ & $=$ & 
$ \{ (2,2)_{17}, (2,2)_{18}, (2,2)_{19}, (2,2)_{20}, (2,2)^N_{46}, {\color{red}(2,2)^N_{47}},  (2,2)^N_{48}, (2,2)^N_{49}, (2,2)^N_{50}, 
$\\ 
&& $  (2,2)_{70}, (2,2)^N_{71} , (2,2)_{ 72} \}$. \\ 
$\overline{ (2,2)^G_{28} } $ & $=$ & 
$ \{ (2,2)_{23}, {\color{red}(2,2)_{24}}, (2,2)_{25}, {\color{red}(2,2)_{26}}, (2,2)_{27}, (2,2)_{28},
(2,2)^N_{46},
{\color{red} (2,2))^N_{47}}, (2,2)^N_{48}, 
$\\ 
&& $ (2,2)^N_{49},  (2,2)^N_{50}, (2,2)_{ 70}, (2,2)^N_{71} ,  (2,2)^N_{72}  \} $. \\ 
$\overline{ (2,2)^G_{38} } $ & $=$ & 
$ \{ 
(2,2)_{35}, (2,2)_{36}, (2,2)_{37}, (2,2)_{38}, (2,2)^N_{46},
 {\color{red} (2,2)^N_{47}},(2,2)^N_{48},
(2,2)^N_{49}, (2,2)^N_{50} 
$,\\ 
&& $(2,2)^N_{68}, (2,2)^N_{69}, (2,2)_{70},(2,2)^N_{71}, (2,2)^N_{72} \} $. \\ 
$\overline{ (2,2)^G_{45} } $ & $=$ & 
$ \{ 
{\color{red}(2,2)_{32}},
(2,2)_{40}, (2,2)_{42},
{\color{red} (2,2)_{41}}, {\color{red}(2,2)_{ 43}}, (2,2)_{44},
(2,2)_{45}, (2,2)^N_{46}, {\color{red} (2,2)^N_{47}}, 
$\\ 
&& $(2,2)^N_{48}, (2,2)^N_{50}, (2,2)^N_{69}, (2,2)_{ 49},  (2,2)^N_{68}, (2,2)_{70}, (2,2)^N_{71}, (2,2)^N_{72} 
\} $. \\ 
$\overline{ (2,2)^G_{52} } $ & $=$ & 
$ \{ 
(2,2)_{ 52}, (2,2)^N_{72}
\} $. \\ 
$\overline{ (2,2)^G_{56} } $ & $=$ & 
$ \{(2,2)_{53}, (2,2)_{54}, (2,2)_{55},
(2,2)_{ 56}, (2,2)^N_{68}, (2,2)^N_{69}, (2,2)_{70}, (2,2)^N_{71}, (2,2)^N_{72}
\} $.  \\ 
$\overline{ (2,2)^G_{57} } $ & $=$ & 
$ \{ 
(2,2)_{57}, (2,2)^N_{68}, (2,2)^N_{72}
\} $. \\ 
$\overline{ (2,2)^G_{59} } $ & $=$ & 
$ \{ 
(2,2)_{53}, (2,2)_{54}, (2,2)_{58}, 
(2,2)_{59}, (2,2)^N_{68}, (2,2)^N_{69}, (2,2)_{70}, (2,2)^N_{71}, (2,2)^N_{72}
\} $.  \\ 
$\overline{ (2,2)^G_{60} } $ & $=$ & 
$ \{ 
(2,2)_{60}, (2,2)^N_{72}
\} $. \\ 
$\overline{ (2,2)^G_{64} } $ & $=$ & 
$ \{ 
(2,2)_{61},
(2,2)_{62}, (2,2)_{63}, (2,2)_{64}, (2,2)^N_{68}, (2,2)^N_{69},  (2,2)_{70}, (2,2)^N_{71},  (2,2)^N_{72}
\} $.  \\ 
$\overline{ (2,2)^G_{66} } $ & $=$ & 
$ \{ (2,2)_{61}, (2,2)_{62}, (2,2)_{65},
(2,2)_{66}, (2,2)^N_{68}, (2,2)^N_{69}, (2,2)_{70}, (2,2)_{ 71}, (2,2)_{ 72}
\} $.  \\ 
$\overline{ (2,2)^G_{67} } $ & $=$ & 
$ \{ 
(2,2)_{67}, (2,2)^N_{72} 
\} $. \\ 
\end{longtable}

\end{theorem}

The irreducible components of $\mathcal{JS}^{(2,2)}$ are represented in Figure \ref{figure:graphic_(2,2)_NA}. In Hasse diagrams we adopt the following notation: the blue color indicates an associative superalgebra, and a square represents a nilpotent superalgebra. Furthermore, for abbreviation, let $i$ stand for either $(2,2)_i$ or $(2,2)_i^N$.

\begin{figure}
 \centering
 \includegraphics[width=20cm, height=10cm, angle=90]{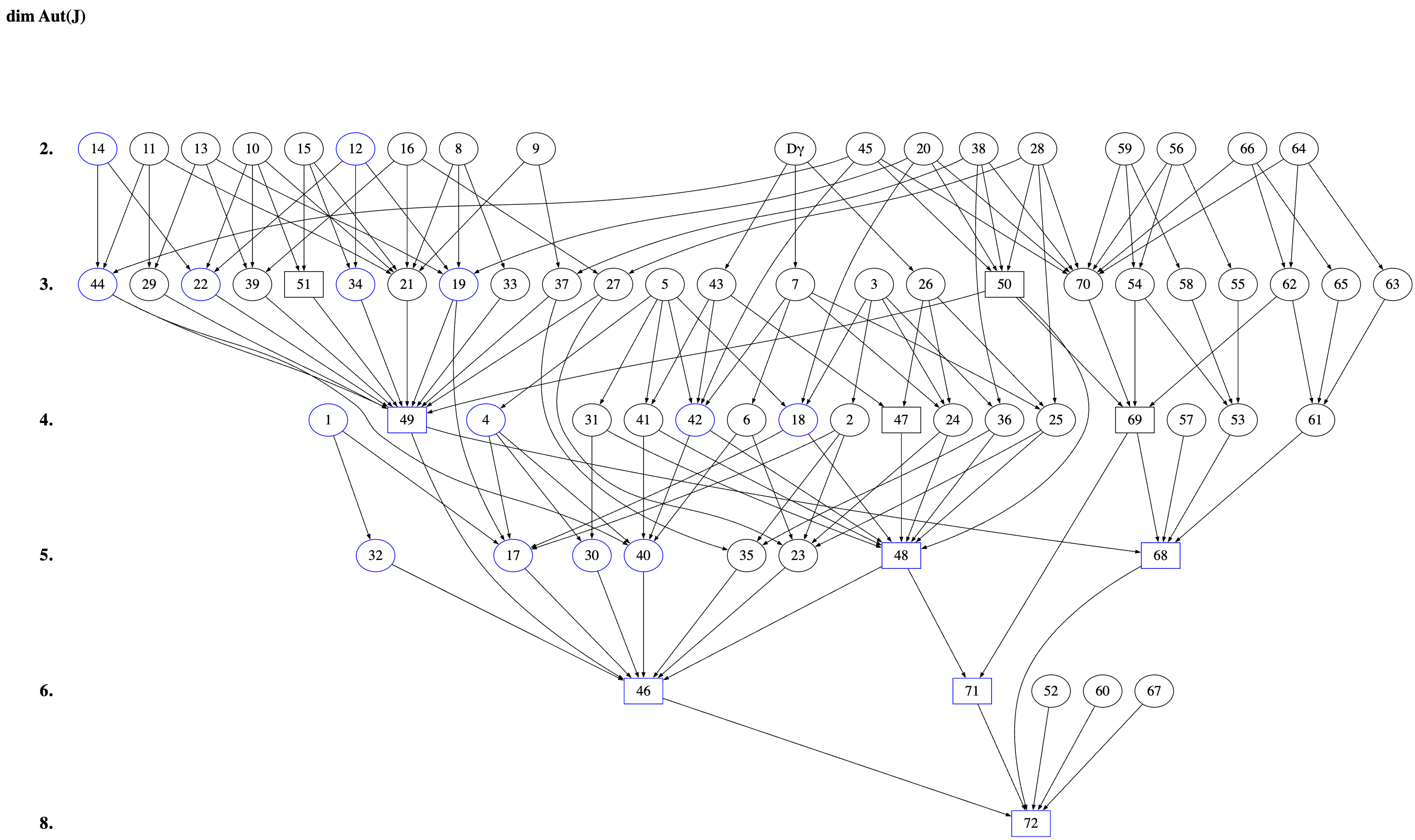}
 \caption{Graphic $(2|2)$}
 \label{figure:graphic_(2,2)_NA}
\end{figure}


Finally, as a direct consequence of Theorem \ref{Thm: principal}, we obtain the following results:

\begin{corollary}
The subvariety $\mathcal{ASC}^{(2,2)} \subset \mathcal{JS}^{(2,2)}$ of supercommutative associative superalgebras of type $(2,2)$ has $6$ irreducible components given by 
\begin{align*}
\overline{ (2,2)_{1}^G }  = & \{ (2,2)_{1}, (2,2)_{17}, (2,2)_{32}, (2,2)^N_{46}, (2,2)^N_{72} 
\}.\\
\overline{ (2,2)_{4}^G } = & \{ 
(2,2)_{4}, (2,2)_{17}, (2,2)_{30}, (2,2)_{40}, (2,2)^N_{46}, (2,2)^N_{72}\} .\\
\overline{ (2,2)_{12}^G } = & \{ 
(2,2)_{12}, (2,2)_{17}, (2,2)_{19}, (2,2)_{22}, (2,2)_{34}, (2,2)^N_{46}, (2,2)^N_{49}, (2,2)^N_{68}, (2,2)^N_{72} \} .\\
\overline{ (2,2)_{14}^G } = & \{ 
(2,2)_{14}, (2,2)_{22}, (2,2)_{40}, (2,2)_{44}, (2,2)^N_{46}, (2,2)^N_{49}, (2,2)^N_{68}, (2,2)^N_{72} \} .\\
\overline{ (2,2)_{18}^G } = & \{ 
(2,2)_{17}, (2,2)_{18}, (2,2)^N_{46},  (2,2)^N_{48}, (2,2)^N_{71}, (2,2)^N_{72}  \}.\\
\overline{ (2,2)_{42}^G } = & \{ 
(2,2)_{40},  (2,2)_{42}, (2,2)^N_{46}, (2,2)^N_{48}, (2,2)^N_{71}, (2,2)^N_{72}
\} .\\
\end{align*}
\end{corollary}

\begin{corollary}
The subvariety $\mathcal{NJS}^{(2,2)} \subset \mathcal{JS}^{(2,2)}$ of nilpotent Jordan superalgebras of type $(2,2)$ has $3$ irreducible components given by
\begin{align*}
\overline{ ((2, 2)^N_{47})^G } = & \{ 
(2, 2)^N_{46}, (2, 2)^N_{47}, (2, 2)^N_{48}, (2, 2)^N_{71}, (2, 2)^N_{72}
\}.\\
\overline{ ((2, 2)^N_{50})^G } = & \{ 
(2, 2)^N_{46}, (2, 2)^N_{48},
(2, 2)^N_{49}, 
(2, 2)^N_{50}, (2, 2)^N_{68}, (2, 2)^N_{69}, (2, 2)^N_{71}, (2, 2)^N_{72}
\}.\\
\overline{ ((2, 2)^N_{51})^G } = & \{ 
(2, 2)^N_{46}, (2, 2)^N_{49}, (2, 2)^N_{51}, (2, 2)^N_{68}, (2, 2)^N_{72}
\}.
\end{align*}
\label{corollary:irreducible_components_nilpotent_3_1}
\end{corollary}

We finish the paper with the open problems in the case (2,2). 

\begin{table}[ht]
\begin{center}
\label{Open_problems_2_2}
\begin{spacing}{1.4}
\begin{tabular}{| l | c |}
\hline
$ 3 \rightarrow 25, 47$; 
$ \;\; 7 \rightarrow 41, 47$; 
$\;\; 8 \rightarrow 51 $; &
\\
$i \rightarrow 47$, for $i \in \{5, 20, 38 \} $; 
$\;\; i \rightarrow 32$, for $i \in \{10, 14 \} $; &
\\
$ i \rightarrow \{ 33, 51 \}$, for $i \in \{11, 13 \} $; $\;\; 12 \rightarrow  32 $;  
& Open Problems
\\
$28 \rightarrow \{ 24, 26, 47 \}$;
$\;\;45 \rightarrow \{32, 41, 43, 47 \}$&
\\
$D_{\gamma}\to i$, for $i\in \{2,4,17,18,30,31,32,35,36\}$ &\\
\hline
\end{tabular}
\end{spacing}
\caption{Open Problems in the case $(2, 2)$}
\end{center}
\end{table}

\bibliographystyle{amsplain}
\bibliography{library}

\end{document}